\documentclass[3p, number,sort&compress,times,fleqn]{elsarticle}      

\usepackage{./csml}

\usepackage{braket}

\usepackage[colorlinks = true,
linkcolor = blue,
urlcolor  = blue,
citecolor = blue,
anchorcolor = blue]{hyperref}

\usepackage{tikz}
\usetikzlibrary{quantikz2}
\usetikzlibrary{calc,positioning,arrows.meta,fpu}
\tikzset{
  step arrow/.style={
    -{Latex[length=2.5mm,width=1mm]},
    black,
    line width=0.4pt
  },
    separator/.style={
    line width=0.2pt
  },
  step label/.style={
    above=2pt,
font=\footnotesize
  }
}
\usepackage[most]{tcolorbox}
\newtcolorbox[auto counter]{eqbox}[2][]{
  enhanced,
  float,
  colframe=csmlDarkBlue!7.5,
  colback=csmlDarkBlue!7.5,
  colbacktitle=csmlDarkBlue!25,
  coltitle=black, 
  boxrule=0.4pt, arc=2pt,
  left=4pt,right=4pt,top=4pt,bottom=4pt,
  boxsep=3pt, % inner padding
  before upper={
    \setlength{\parskip}{0pt}%
    \setlength{\abovedisplayskip}{6pt}%
    \setlength{\belowdisplayskip}{6pt}%
    \setlength{\abovedisplayshortskip}{4pt}%
    \setlength{\belowdisplayshortskip}{4pt}%
    \setlength{\jot}{2pt}% aligns' interline space
  },
  title={Box~\thesection.\thetcbcounter: #2},
  #1
}
\renewcommand{\coloneq}{\mathrel{\vcentcolon=}}

\usepackage{booktabs}

\graphicspath{{./figs/}}
\begin{document}
	
\begin{frontmatter}

\title{QAFE$^2$: Quantum accelerated multiscale finite element analysis}

\author[1]{Yiren Wang}
\author[2,3]{Michael Ortiz}
\author[1]{Fehmi Cirak\corref{cor1}}
\ead{f.cirak@eng.cam.ac.uk}

\cortext[cor1]{Corresponding author}

\address[1]{Department of Engineering, University of Cambridge, Cambridge, CB2 1PZ, UK }
\address[2]{Division of Engineering and Applied Science, California Institute of Technology, Pasadena, CA 91125, USA}
\address[3]{Centre Internacional de M\`etodes Num\`erics a l'Enginyeria (CIMNE), Universitat Polit\`ecnica de Catalunya, 08034 Barcelona, Spain}

\begin{abstract}
The computational cost of concurrent multiscale finite element methods is dominated by the repeated solution of microscopic representative volume element (RVE) problems at macroscopic quadrature points. In this work, we introduce a quantum--classical framework for multiscale finite element analysis (\emph{QAFE$^2$}) that leverages quantum parallelism to fundamentally alter the scaling of RVE-based homogenisation. At the single-RVE level, the proposed quantum solver attains polylogarithmic complexity with respect to the microscopic discretisation size, yielding an exponential asymptotic speedup over the best available classical solvers. More importantly, \emph{QAFE$^2$} exploits quantum superposition and entanglement to evaluate, in a single quantum execution, the entire ensemble of RVE problems associated with all macroscopic quadrature points. This capability is a form of intrinsic quantum concurrency with no classical analogue. Numerical experiments on one- and two-dimensional model problems with known analytical solutions confirm the accuracy of the proposed formulation and verify the theoretical computational scaling and parallel performance. 
\end{abstract}


\end{frontmatter}

%\tableofcontents
%
%--------------------------------------------------------------------------------
\section{Introduction \label{sec:intro}}
%--------------------------------------------------------------------------------
%

% Introduction, motivation, and scope

The progress of computational mechanics over the past half-century owes largely to the sustained growth of classical computing power. From its origins in the late 1950s, the field has advanced by exploiting increases in problem size made possible by faster processors and greater memory, culminating in finite element simulations of unprecedented scale. This trajectory, however, is now approaching fundamental physical limits. The slowing of Moore’s law, coupled with escalating energy and infrastructure costs, makes clear that further gains in predictive capability cannot rely on brute-force classical scaling alone. A qualitative change in computing paradigms is therefore required.

Quantum computing offers precisely such a paradigm shift. By exploiting superposition, entanglement, and unitary time evolution, quantum computers process information in a manner that is fundamentally distinct from classical digital computers. Of particular significance is the ability of a quantum system to evolve an entire state vector simultaneously, giving rise to quantum parallelism and interference. The central question addressed in this work is whether these distinctive features can be harnessed in a meaningful way to accelerate core tasks in computational mechanics, and in particular those arising in concurrent multiscale analysis.
\begin{figure}
\centering
\includegraphics[width=\textwidth]{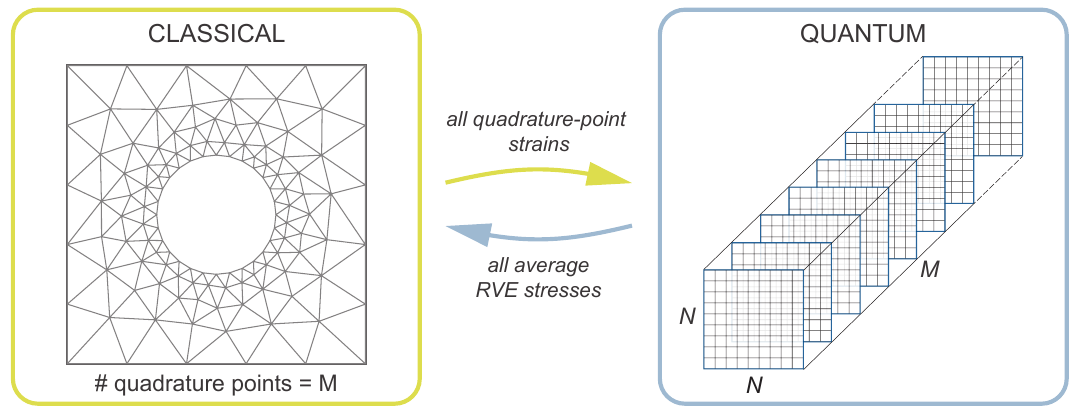} 
\caption{Schematic of the \emph{QAFE$^2$} framework for multiscale analysis.  The macroscale problem is solved on a classical computer, while the microscale RVE  problems associated with the~$M$ quadrature points are solved on a quantum computer.  In \emph{QAFE$^2$}, the $M$ microscale problems are each discretised with an $N \times N$ uniform grid and are solved simultaneously in a single quantum computation. \label{fig:intro_qafe2}}
\end{figure}

% Historical Background

The origins of quantum computing can be traced to the early 1980s. Feynman argued that classical computers are fundamentally ill-suited for simulating quantum systems, since the dimensionality of the quantum state space grows exponentially with the number of particles~\cite{feynman1982}. Closely related arguments were advanced independently by Manin~\cite{manin1980}. From a complementary perspective, Benioff emphasised the thermodynamic cost of irreversible classical logic and proposed computation based on reversible quantum dynamics~\cite{benioff1982}. These ideas were unified and formalised by Deutsch, who introduced the concept of a universal quantum computer and established the theoretical basis for quantum algorithms that outperform their classical counterparts~\cite{deutsch1985}.

Subsequent breakthroughs, notably Shor’s polynomial-time algorithm for integer factorisation~\cite{shor1997} and Grover’s quadratic-speedup search algorithm~\cite{grover1996}, demonstrated unequivocally that quantum advantage is not limited to the simulation of quantum physics. Nevertheless, the translation of these insights into computational mechanics has been limited. Much of the existing literature focuses on quantum linear system solvers, particularly the seminal Harrow–Hassidim–Lloyd (HHL) algorithm~\cite{harrow2009}. Later linear system solvers, offering improved scaling and resource efficiency, include the approaches of Childs et al.~\cite{childs2017} and Gily\'en et al.~\cite{gilyen2019quantum}; see also the recent review~\cite{morales2024quantum}. While theoretically appealing, such linear system solvers rely on stringent assumptions regarding sparsity, conditioning, and efficient amplitude encoding of data, conditions that are rarely met in large-scale mechanics applications~\cite{aaronson2015read,montanaro2016quantum}. This gap motivates the search for alternative formulations that more directly exploit the structure of mechanics problems.

% Multiscale Modeling and the RVE Bottleneck

The present work focuses specifically on multiscale modelling of materials, where macroscopic structural behaviour emerges from the collective response of complex microstructures~\cite{phillips2001crystals,ortiz2001mixed,bhattacharya2003microstructure}. In many concurrent multiscale schemes---such as computational homogenisation and FE$^2$ methods~\cite{terada1997, kouznetsova2001, miehe2002}---the constitutive response at each macroscopic material point is obtained by solving a boundary value problem on a microscopic RVE~\cite{suquet1987, geers2010}. While such approaches are theoretically sound and physically transparent, they are computationally intractable: the repeated, on-the-fly solution of RVE problems typically dominates the overall cost.

FFT-based homogenisation methods, pioneered by Moulinec and Suquet~\cite{moulinec1998} and subsequently refined by many authors~\cite{eyre1999, willot2015}, have significantly reduced the cost of periodic RVE solutions, achieving near linear complexity~$O(N\log N)$; see also the reviews~\cite{schneider2021review, lucarini2022fft}. Despite these advances, the sheer number of RVEs required in fully concurrent simulations renders even FFT-based schemes computationally prohibitive at scale. This persistent bottleneck has effectively limited the routine application of FE$^2$ methods in large-scale engineering simulations.

% Objectives and Main Results

The principal objective of the present work is to demonstrate that quantum computing can fundamentally change the computational scaling of RVE-based multiscale analysis. We specifically build on the work of Liu \emph{et al.} \cite{liu2024towards}, who developed a quantum reformulation of the classical FFT-based fixed-point scheme of Moulinec and Suquet~\cite{moulinec1998} for one-dimensional RVE problems in antiplane shear. In this setting, the unknown field is scalar-valued, and the Fourier-space Green's operator is a diagonal matrix. The key observation of that work is that the Quantum Fourier Transform (QFT) performs the same mathematical operation as the classical FFT, but with polylogarithmic complexity $O(\log^2 N)$ \cite{coppersmith2002,ikeAndMike}. Building on this foundation, Liu~\emph{et al.} \cite{liu2024towards} introduced quantum circuits for polynomial and piecewise-Chebyshev representations of material heterogeneity~\cite{woerner2019quantum, vazquez2022enhancing}, algebraic operations in Fourier space, and a measurement-efficient implementation of fixed-point iterations. When assembled into a complete solver, these components yield an overall complexity $O(\log^c N)$ for some modest constant $c$, representing an exponential speedup relative to the best available classical algorithms. Theoretical complexity estimates are supported by numerical experiments that confirm the anticipated scaling behaviour. Extensions to non-periodic boundary value problems have been proposed in recent work~\cite{febrianto2025quantum}. 

Building on those foundations, we propose a quantum--classical framework for concurrent multiscale finite element analysis, termed \emph{QAFE$^2$}, for multidimensional problems requiring the repeated solution of RVE problems. Hence, \emph{QAFE$^2$}  targets the dominant computational bottleneck of FE\textsuperscript{2}-type methods: the solution of RVE problems at macroscopic quadrature points. The approach reformulates the classical FFT-based homogenisation scheme of Moulinec and Suquet into a fully quantum algorithm by exploiting band-limited Fourier representations, fixed-point iterations in Fourier space, and the QFT.  In multidimensional problems, the unknown field is vector-valued, and the Fourier-space Green's operator is a block-diagonal matrix. Central to \emph{QAFE$^2$} is a block-encoded implementation of the Fourier-space Green's operator via the linear combination of unitaries (LCU) technique~\cite{childs2012hamiltonian,herbert2026quantum}, enabling non-unitary constitutive updates to be realised within unitary quantum circuits. In contrast to quantum linear-system solvers that rely on restrictive assumptions, \emph{QAFE$^2$} mirrors the structure of established computational mechanics algorithms while leveraging quantum parallelism in a direct and physically transparent manner.

A chief accomplishment of the present work is the demonstration that \emph{QAFE$^2$} affords an exponential reduction in asymptotic complexity for RVE solves relative to state-of-the-art classical FFT-based schemes. For an RVE discretised on an $N\times N$ grid, the single-RVE quantum algorithm achieves polylogarithmic complexity in $N$, replacing the classical $\mathcal{O}(N^2\log N)$ scaling by $\mathcal{O}(\log^c N)$ for a modest constant $c$. Beyond single-RVE acceleration, the \emph{QAFE$^2$} solver sets forth a fully concurrent quantum treatment of all $M$ RVEs associated with a macroscopic finite element mesh, see Figure~\ref{fig:intro_qafe2}. By encoding the ensemble of macroscopic strains into orthogonal subspaces of a single quantum state, all $M$ RVEs are solved simultaneously using a single instance of the fixed-point iteration circuit. The resulting overall complexity scales as $\mathcal{O}(M\,\log^c M + \log^c N)$, representing a qualitative departure from classical FE\textsuperscript{2}, where the cost grows linearly with $M$ times the single-RVE solve, see Table~\ref{tab:complexity_scaling}.
\begin{table}
\centering
\caption{Computational complexity of classical and quantum FE$^2$ multiscale frameworks, considering $M$ macroscopic quadrature points and a microscopic RVE discretisation of $N \times N$ grid points. The estimates assume that a separate RVE solve is required at each quadrature point. The exponent $c$ is an implementation-dependent constant that may differ across expressions.}
\label{tab:complexity_scaling}
\addtolength{\tabcolsep}{6pt} % Improves readability by spacing out columns
\begin{tabular}{@{}lccc@{}}
\toprule
 & {Initialisation} & {Computation} & {Total } \\ 
\midrule
Classical               & $O(1)$                  & $O(M N^{c})$               & $O(M N^{c})$                    \\
Sequential quantum     & $O(1)$                  & $O(M \log^{c} N)$        & $O(M \log^{c} N)$             \\
Parallel quantum (\emph{QAFE$^2$})        & $O(M\log^{c} M)$             & $O(\log^{c} N)$        & $O(M\log^{c} M + \log^{c} N)$    \\ 
\bottomrule
\end{tabular}
\end{table}

The theoretical analysis is supported by numerical experiments on one- and two-dimensional model problems in antiplane shear with known analytical solutions. These examples confirm the predicted scaling behaviour, demonstrate the accuracy of the quantum fixed-point iteration, and illustrate the effectiveness of quantum parallelisation across multiple macroscopic loading states. Although only linear model problems are considered and the simulations are carried out on noiseless quantum emulators, they provide concrete evidence that \emph{QAFE$^2$} admits explicit circuit constructions with gate counts consistent with the theoretical estimates.

The paper is organised as follows. Section~\ref{sec:model_problem} introduces the model problem of antiplane shear and reviews the classical FFT-based fixed-point formulation, together with its band-limited Fourier discretisation. Section~\ref{sec:quantum_serial} details the quantum implementation of the two-dimensional model problem, including the Fourier-space solution, incremental strain updates, and the construction of the fixed-point iteration circuit. Section~\ref{sec:quantum_parallel} extends the formulation to the simultaneous quantum solution of all RVEs arising in a macroscopic finite element computation and analyses the resulting parallel complexity. Section~\ref{sec:examples} presents numerical examples assessing accuracy and scaling, and Section~\ref{sec:conclusions} concludes with a discussion of implications, limitations, and directions for future work in quantum-accelerated computational mechanics. Additional details are provided in two appendices.

%
%--------------------------------------------------------------------------------
\section{Model problem: Antiplane shear \label{sec:model_problem}}
%--------------------------------------------------------------------------------
%
In multiscale finite element analysis, the constitutive response at each quadrature point of the macroscopic FE discretisation is determined by solving an associated microscopic RVE problem. The macroscopic strain is imposed as the average strain of the RVE, and the average stress of the RVE represents the homogenised macroscopic stress.  For the purposes of exposition, we consider an inhomogeneous two-dimensional elastic solid undergoing antiplane shear. The proposed scheme applies broadly to general RVE problems.  

%--------------------------------------------------------------------------------
\subsection{Problem formulation \label{sec:problem_formulation}}
%--------------------------------------------------------------------------------
%
We  consider the RVE domain~\mbox{$\Omega = (0, \, L) \times (0, \, L) \in \mathbb R^2 $} with an edge length of~$L$. The coordinates of the points \mbox{$ \vec x \in \Omega$} are denoted as~\mbox{$\vec x = ( x_0 \quad x_1 )^\trans$}. The shear modulus~$\mu (\vec x) \in \mathbb R$ in the RVE is periodic, i.~e., 
\begin{equation}
	\mu(x_0, x_1) = \mu (x_0 + L,  x_1) = \mu (x_0 ,  x_1+L) \, . 
\end{equation}
The RVE is subject to a uniform average strain vector~$\overline { \vec \gamma} \in \mathbb R^2$ passed down by the macroscopic problem. The deformation of the RVE is characterised by a scalar transverse displacement field~$u( \vec x) \in \mathbb R$ with shear-strain vector 
\begin{equation} \label{eq:shearStrain}
	 \vec \gamma (\vec x) = \nabla u (\vec x) \, ,  
\end{equation}
where~$\nabla$ denotes the gradient operator and~$ \vec \gamma (\vec x) \in \mathbb R^2$. The displacement field~$u(\vec x) $ consists of an affine component matching the prescribed average strain~$\overline {\vec  \gamma}$ and a fluctuating component~$v(\vec x) \in \mathbb R$, so that 
\begin{equation} \label{eq:dispDecomposition}
	u(\vec x) = \overline{ \vec \gamma } \cdot \vec x + v( \vec x) \, , 
\end{equation}
which implies the strain decomposition 
\begin{equation} \label{eq:strainDecomposition}
	\vec \gamma (\vec x) = \overline{ \vec \gamma } + \nabla v( \vec x) \, .
\end{equation}
The fluctuating displacement~$v (\vec x)$ must satisfy the periodicity condition
\begin{subequations} \label{eq:dispPeriodicitiy}
\begin{align} 
  v(x_0,  0) = v(x_0,  L) \, , 
  & \quad
  x_0 \in (0, \, L) \, , 
  \\ 
  v(0, x_1) = v(L, x_1) \, ,
  & \quad
  x_1 \in (0, \, L)  \, .
\end{align}
\end{subequations}
Therefore, 
\begin{equation} \label{eq:dispZeroAverage}
 	 \int_\Omega  \nabla v( \vec x) \D \vec x = \vec 0 \, .
\end{equation}
The strain decomposition~\eqref{eq:strainDecomposition} together with the zero-average condition~\eqref{eq:dispZeroAverage} ensure that~$\overline{ \vec \gamma}$ is indeed the average strain. We verify that
\begin{equation}
	\frac{1}{L^2} \int_\Omega \vec \gamma (\vec x)  \D \vec x  = \frac{1}{L^2} \int_\Omega \left ( \overline{ \vec \gamma } + \nabla v( \vec x)  \right ) \D \vec x  = \overline{\vec \gamma} \, .
\end{equation}
The stress field vector~$ \vec \sigma (\vec x) \in \mathbb R^2$ of the RVE must satisfy the equilibrium equation 
\begin{align} \label{eq:equilibrium}
	\nabla \cdot  \vec \sigma (\vec x) = 0 \, , 	
\end{align}
and the constitutive equation 
\begin{equation} \label{eq:constitutive}
	\vec \sigma( \vec x) = \mu (\vec x)  \vec \gamma (\vec x) = \mu (\vec x)  \left ( \overline{ \vec \gamma } + \nabla v( \vec x) \right ) \, .
\end{equation}
Hence, the boundary value problem for the microscopic RVE can be summarised as 
\begin{equation} \label{eq:homogBVP}
	\nabla \cdot \left ( \mu ( \vec x)   \nabla v( \vec x) \right ) + \overline { \vec \gamma} \cdot \nabla \mu (\vec x) = 0  \, ,
\end{equation}
subject to the periodicity condition~\eqref{eq:dispPeriodicitiy}. 

Due to the non-constant shear modulus~$\mu( \vec x)$, it is not possible to solve the boundary value problem~\eqref{eq:homogBVP} directly using the Fourier transform. Therefore, following Moulinec and Suquet~\cite{moulinec1998} we introduce a constant reference shear modulus~$\mu_0$ and rewrite the constitutive equation~\eqref{eq:constitutive} as 
\begin{equation} \label{eq:homogBVPlin}
    \vec \sigma( \vec x)  = \left ( \mu ( \vec x) - \mu_0 \right )  \vec \gamma( \vec x) + \mu_0 \vec  \gamma( \vec x)   =    \vec \tau (\vec x) + \mu_0 \vec \gamma(\vec x)	\, , 
\end{equation}
where~$\vec \tau( \vec x)$ is referred to as the polarisation stress. Next, we use the equilibrium condition~\eqref{eq:equilibrium} to formulate the fixed-point iteration 
\begin{equation} \label{eq:homogBVP_it}
	 \mu_0 \nabla \cdot \nabla v^{(s+1)}(\vec x) + \nabla \cdot  \vec \tau^{(s)} ( \vec x) = 0 \, . 
\end{equation}
The polarisation stress~$\vec \tau^{(s)}$ at iteration step~$(s)$ depends on the known displacement~$v^{(s)}(\vec x)$ and the applied strain~$\overline{ \vec \gamma}$.

Finally, the homogenised stress~$\overline{ \vec \sigma}$ for the macroscopic problem is given by 
\begin{equation} 
	\overline{ \vec \sigma} =  \frac{1}{L^2} \int_\Omega  \vec \sigma (\vec x) \D \vec x   =  \frac{1}{L^2} \int_\Omega  \mu(\vec x)  \left (\overline { \vec \gamma} + \nabla v(\vec x) \right )  \D \vec x \, .
\end{equation}

%

%--------------------------------------------------------------------------------
\subsection{Band-limited Fourier discretisation  \label{sec:band_limited_approximation}}
\begin{figure}
\centering
	\includegraphics[scale=0.95]{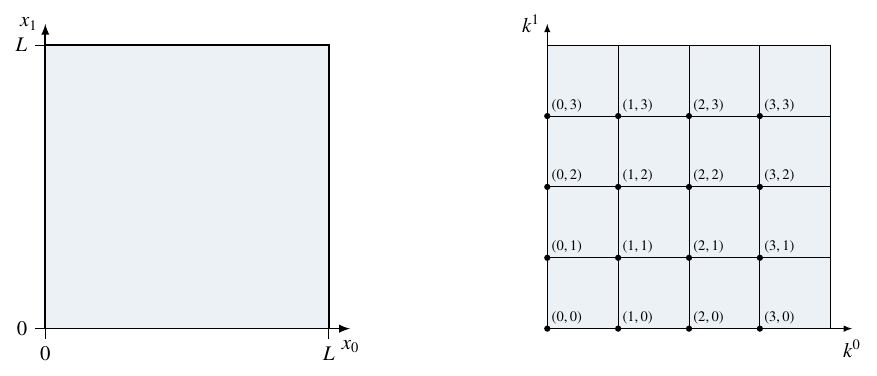} 
	 \caption{RVE domain~$\Omega = (0, \, L) \times (0, \, L)$ and its discretisation with~$N \times N$ cells, where $N=4$. The grid points are labelled as~$\vec {k }= (k^0, \, k^1)$. The problem is periodic in both directions. \label{fig:2D_grid_labelling}}
\end{figure}

%--------------------------------------------------------------------------------
%
We represent the periodic fluctuation field~$v(\vec x)$ over the RVE domain $\Omega=(0, \, L)\times(0, \, L)$ using the band-limited Fourier approximation 
\begin{equation} \label{eq:fourierApprox-2d}
	v (\vec x) \approx v^h(\vec x) =  \frac{1}{N} \sum_{k^0, \, k^1=-N/2}^{N/2-1} \hat v_{\vec k} e^{\frac{2 i \pi}{L} \vec k \cdot \vec x  }  \, , 
\end{equation}
where $N$ is an even positive integer and  $\vec k = ( k^0, \, k^1)$ is a multi-index and~\mbox{$\vec k \cdot \vec x \equiv k^0 x_0 + k^1 x_1$}. We rewrite this as 
\begin{equation} \label{eq:fourierApprox}
	v^h(\vec x) =  \frac{1}{N} \sum_{k^0, \, k^1=0}^{N-1}  \hat v_{\vec k} e^{i \vec \xi_{\vec k} \cdot \vec x }  \, , 
\end{equation}
with the discrete wave vectors 
\begin{equation}
	\vec \xi_{\vec k} 
  = 
  \begin{pmatrix}
    \dfrac{2 \pi r(k^0)}{L} & 
    \dfrac{2 \pi r(k^1)}{L} 
  \end{pmatrix}^\trans \, ,
\end{equation}
and the relabelling function
\begin{equation} \label{eq:componentsToFrequencies}
	r(k) = \begin{cases}
		k & 0 \le k < N/2 \\
		k - N & N/2 \le k < N \, .
	\end{cases}
\end{equation}
The coefficients~$\hat v_{\vec k} $ correspond to the discrete Fourier transform (DFT) of the grid-point values~$v^h (\vec x_{\vec k})$ sampled over the uniform grid
\begin{equation}
	\vec x_{\vec k} = \frac{L}{N}  \begin{pmatrix} k^0 & k^1 \end{pmatrix}^\trans \, , \quad  k^0,k^1 \in \{0, \, 1, \, \dotsc, \, N-1\} \, ; 
\end{equation}
see Figure~\ref{fig:2D_grid_labelling}. 

After approximating the periodic fields within the fixed-point iteration~\eqref{eq:homogBVP_it} using the band-limited expansion~\eqref{eq:fourierApprox} and invoking the orthogonality of the Fourier basis, the incremental RVE  solution is obtained as
\begin{equation}
	v^{h, (s+1)} (\vec x) = \frac{1}{\mu_0}  \sum_{k^0, \, k^1=0}^{N-1}  \frac{ i \vec \xi_{\vec k} \cdot \hat{\vec \tau}_{\vec k}^{(s)}}{\vec \xi_{\vec k} \cdot \vec \xi_{\vec k}} e^{i \vec \xi_{\vec k} \cdot \vec x} \, ,\qquad \forall {\vec \xi_{\vec k}} \neq \vec 0 \, .
\end{equation}
Finally, after differentiation, the fluctuation strain is given by 
\begin{equation} \label{eq:fixPointIterStrain}
	\nabla  v^{h , (s+1)} (\vec x) = -  \frac{1}{\mu_0}  \sum_{k^0, \, k^1=0}^{N-1}  \frac{ \vec \xi_{\vec k}  \vec \xi_{\vec k}^\trans }{\vec \xi_{\vec k} \cdot \vec \xi_{\vec k}}  \hat{\vec \tau}_{\vec k}^{(s)}   e^{i \vec \xi_{\vec k} \cdot \vec x}
	=   \sum_{k^0, \, k^1=0}^{N-1} \hat {\vec \Gamma}_{\vec k}  \hat {\vec \tau}_{\vec k}^{(s)}  e^{i \vec \xi_{\vec k} \cdot \vec x} \, ,
	 \qquad \forall {\vec \xi_{\vec k}} \neq \vec 0  \, .
\end{equation}
We refer to~$ \hat {\vec \Gamma}_{\vec k} \in \mathbb R^{2\times 2}$ in the following as the (Fourier-space) strain Green's operator.

%
%--------------------------------------------------------------------------------
\section{Quantum implementation \label{sec:quantum_serial}}
%--------------------------------------------------------------------------------
%
We proceed to the quantum implementation of the band-limited Fourier discretisation introduced in the preceding section. In this section, we switch to Dirac notation to distinguish between quantum mechanical and classical vectors. Introductions to Dirac notation can be found in standard texts on quantum computing~\cite{ikeAndMike,herbert2026quantum} and earlier work~\cite{liu2024towards}.  
 %
%--------------------------------------------------------------------------------
\subsection{Overview \label{sec:implementation_overview}}
%--------------------------------------------------------------------------------
%
Recall that the microscopic RVE problems take as input the macroscopic FE strain~\mbox{$\ket{\overline \gamma} = \overline {\gamma}_0 \ket 0  +  \overline {\gamma}_1 \ket 1$} and return the average stress~\mbox{$\ket{\overline \sigma}  = \overline {\sigma}_0 \ket 0  + \overline {\sigma}_1\ket 1$}. All RVEs are discretised using a uniform grid with~$N \times N$ cells, where the grid points are indexed by the multi-index $\vec {k}=(k^0, \,  k^1) $. The grid-point strain vector~$\ket \gamma$, shear modulus vector~$\ket \mu$ and polarisation stress~$\ket \tau$ are expressed as
\begin{subequations}
\begin{align}
	\ket \gamma &= \sum_{k^0, \, k^1 =0}^{N-1} \ket{k^0}\ket{ k^1}\ket {\gamma_{k^0k^1}}   = \sum_{k^0, \, k^1=0}^{N-1} \left( \gamma_{0, \, k^0k^1} \ket{k^0}\ket{k^1}\ket 0 +  \gamma_{1, \, k^0k^1} \ket{k^0}\ket{k^1}\ket 1 \right)  \, ,\\
	\ket \mu &= \sum_{k^0, \, k^1=0}^{N-1}  \mu_{k^0k^1} \ket{k^0}\ket{k^1}  \, , \\ 
	\ket \tau &= \sum_{k^0, \, k^1=0}^{N-1} \ket{k^0}\ket{k^1}\ket {\tau_{k^0k^1}}   = \sum_{k^0, \, k^1=0}^{N-1} \left( \tau_{0, \, k^0k^1} \ket{k^0}\ket{k^1}\ket 0 +  \tau_{1, \, k^0k^1} \ket{k^0}\ket{k^1}\ket 1 \right)  \, .
\end{align}
\end{subequations}
Hence,~$\ket \gamma, \ket \tau \in \mathbb R^{2N^2}$  and~$\ket \mu \in \mathbb R^{N^2}$. The strain vector~$\ket \gamma$ is determined by solving the RVE problem via the fixed-point iteration~\eqref{eq:fixPointIterStrain}. As detailed in Box~\ref{tb:rve_problem}, each iteration step consists of five distinct substeps, which are repeated until the average stress~$\ket{\overline \sigma}$ is converged. The average stress is given by 
\begin{equation}
	\ket{\overline \sigma} = \frac{1}{N^2} \sum_{k^0, \, k^1=0}^{N-1} \mu_{k^0k^1}\ket{\gamma_{k^0k^1}} \, .
\end{equation}
\begin{eqbox}{Incremental RVE update \label{tb:rve_problem}}
Each step $(s)$ consists of the following substeps. 
\begin{enumerate}
	\item[S1.] Computation of the polarisation stress:
	\begin{equation*}
		\ket {\tau}^{(s)} = \sum_{k^0, k^1}  (\mu_{k^0k^1} - \mu_0)  \ket{k^0}\ket{k^1}\ket {\gamma_{k^0k^1}}^{(s)}  \, .
	\end{equation*}
	\item[S2.] Quantum Fourier transform of the polarisation stress:
	\begin{equation*}
		\ket {\hat \tau}^{(s)} = QFT \otimes QFT  \ket {\tau}^{(s)}  \, .
	\end{equation*}
	\item[S3.] Computation of the fluctuation strain:  
	\begin{equation*}
		\ket {\hat \gamma}^{(s+1)} =   \sum_{k^0, \, k^1} \ket{k^0}\ket{k^1} \hat \Gamma_{k^0k^1}\ket {\hat \tau_{k^0k^1}}^{(s)}  \, .
	\end{equation*}
	\item[S4.] Application of the macroscopic FE strain: 
	\begin{equation*}
		\hat{ \gamma}^{{(s+1)}}_{0, \, k^0 k^1} =  \bar{ \gamma}_0  \,  ,  \quad \hat{ \gamma}^{{(s+1)}}_{1, \, k^0 k^1} =  \bar{ \gamma}_1\qquad    \text{for $k^0=k^1 =0$} \, .\\
	\end{equation*}
	\item[S5.] Inverse quantum Fourier transform of the strain:
	\begin{equation*}
		\ket {\gamma}^{(s)} = QFT^{\dagger} \otimes  QFT^{\dagger}  \ket {\hat \gamma}^{(s+1)}  \, .
	\end{equation*} 
\end{enumerate}
\end{eqbox}

In devising the quantum circuit, the most subtle step is S3, which will be discussed in Section~\ref{sec:incrementalFourierSpace}. Subsequently, in Section~\ref{sec:incrementalUpdate} we introduce the  quantum circuit for computing one incremental step, and in Section~\ref{sec:fixedPointIter} the circuit for iterating over several steps.  

%
%--------------------------------------------------------------------------------
\subsection{Fourier space solution \label{sec:incrementalFourierSpace}}
%--------------------------------------------------------------------------------
%
We consider the solution of the incremental problem in step S3 and omit the iteration step index $ (s)$ to avoid clutter. For each discrete wave vector, 
\begin{equation}
	\ket {\xi_{k^0k^1} } = \xi_{0, k^0} \ket 0 + \xi_{1, k^1} \ket 1 =  \frac{2 \pi r(k^0)}{L}  \ket 0 + \frac{2 \pi r(k^1)}{L}  \ket 1  \, , 
\end{equation}
the Fourier space strain~$\ket{\hat {\gamma}_{k^0 k^1}} \in \mathbb R^2 $ is given by
\begin{equation}
	\ket{\hat {\gamma}_{k^0 k^1}} =   \hat \Gamma_{k^0k^1}\ket {\hat \tau_{k^0k^1}} = - \frac{1}{\mu_0} \frac{\ket{\xi_{k^0k^1}}  \bra{\xi_{k^0k^1}}}{\braket{\xi_{k^0k^1 }  \mid \xi_{k^0k^1} }} \ket {\hat \tau_{k^0k^1}} \, , 
\end{equation}
where~$\ket {\hat \tau_{k^0k^1}} \in \mathbb R^2$ is the Fourier space polarisation stress.  Since the matrix~$\hat \Gamma_{k^0k^1} \in \mathbb R^{2 \times 2}$  is non-unitary, it cannot be directly quantum encoded. Therefore, we construct a larger unitary matrix~$U_{\hat \Gamma}$ such that~$\hat \Gamma_{k^0k^1} \in \mathbb R^{2 \times 2}$ appears as its one of the subblocks. This process is referred to as block encoding.  

\begin{figure} 
\centering

\subfloat[][$U_{\text{poly}} (\alpha_l)$ ]{ \label{fig:lcu_basic_a} 
	\includegraphics[scale=0.875]{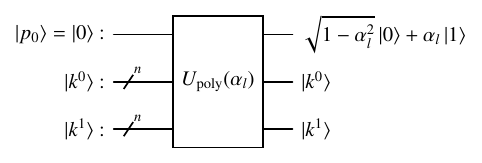}
}
\hspace{0.04\textwidth}
\subfloat[][$U_{\hat \Gamma_l} = U_{\text{poly}} (\alpha_l)  \otimes U_{c_l}$]{  \label{fig:lcu_basic_b}
	\includegraphics[scale=0.875]{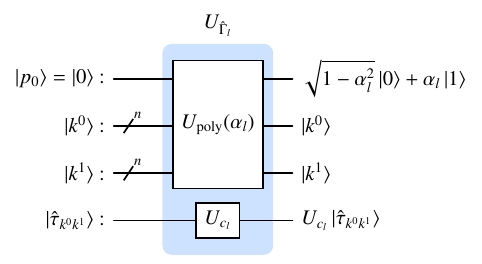}
}
\caption{Quantum circuit for computing the strain~$\ket{\hat \gamma_{k^0k^1}}$ corresponding to one of the terms in the LCU decomposition~\eqref{eq:gammaLCU_0}. The discretisation consists of~$2^n \times 2^n$ grid points. The physically relevant amplitudes of the state vector correspond to basis vectors of the form~$\ket 1 \ket{k^0} \ket{k^1} $ in (a) and~$\ket 1 \ket{k^0} \ket{k^1}  \ket{\hat \tau}$ in~(b).
\label{fig:lcu_basic}}
\end{figure}

One ubiquitous approach for block encoding is the linear combination of unitaries (LCU) technique. To this end, we write the strain Green's operator~$\hat \Gamma_{k^0k^1} $ as a linear combination of the identity matrix~\mbox{$I \in \mathbb R^{2 \times 2}$} and the orthogonal Pauli matrices~\mbox{$X  \in \mathbb R^{2 \times 2}$} and~\mbox{$Z  \in \mathbb R^{2 \times 2}$}, 
\begin{equation} \label{eq:gammaLCU_0}
	\hat \Gamma_{k^0k^1} =   
	 - \frac{1}{2 \mu_0}  I -\frac{1}{\mu_0}  \frac{r(k^0) r(k^1)}{r^2(k^0) + r^2(k^1)}  X - \frac{1}{2\mu_0}   \frac{r^2(k^0) - r^2(k^1)}{r^2(k^0) + r^2(k^1)} Z = 
	 \alpha_0  I + \alpha_1 (k^0, k^1) X  +  \alpha_2 (k^0, k^1) Z = 
	 \sum_{l=0}^2 \alpha_{l} (k^0,  k^1)  U_{c_l}  \, , 
\end{equation}
where~$U_{c_l} \in \{ I, \, X, \, Z\}$. The scalar coefficients~$\alpha_1 (k^0,  k^1)$ and~$\alpha_2 (k^0,  k^1)$ are first approximated classically as polynomials and subsequently quantum encoded following the approach discussed in~\cite{liu2024towards,febrianto2025quantum}; see also~\ref{app:extended-domain}. For simplicity, we retain the notation~$\alpha_l(k^0,k^1)$ for the corresponding polynomial approximants. Each polynomial~$\alpha_l(k^0, k^1)$  is encoded with the help of a suitably constructed unitary~$U_{\text{poly}} (\alpha_l)$ as
\begin{equation} \label{eq:U_poly}
	 U_{\text{poly}} (\alpha_l) \colon \underset{p_0}{\ket 0} \ket {k^0} \ket{ k^1} \mapsto \left ( \sqrt{1-  \alpha_l^2(k^0, k^1)}  \underset{p_0}{\ket{0}} +  \alpha_l (k^0, k^1)    \underset{p_0}{\ket{1}}  \right )    \ket {k^0} \ket{k^1}   \, .
\end{equation}
See Figure~\ref{fig:lcu_basic_a} for the quantum circuit representation of this mapping. Here and in the following, some qubit labels are given in underset notation to improve readability.  The qubit~$\ket{p_0}$, which is initially in state~$\ket 0$, is an ancilla. In the output the sought value~$\alpha_l(k^0,k^1)$  is the amplitude of the state~$\ket 1 \ket{k^0} \ket{k^1}$.  Furthermore, we define \mbox{$U_{\hat \Gamma_l} \coloneq   U_{\text{poly}} (\alpha_l)  \otimes U_{c_l}$}  corresponding to the components of the LCU decomposition~\eqref{eq:gammaLCU_0}. Its action reads  
\begin{equation} \label{eq:U_hat_Gamma}
	U_{\hat \Gamma_l} \colon \underset{p_0}{\ket 0} \ket {k^0}\ket{k^1}  \ket{\hat \tau_{k^0k^1}} \mapsto  \left ( \sqrt{1-  \alpha_l^2(k^0, k^1) }  \underset{p_0}{\ket{0}} +  \alpha_l (k^0, k^1)    \underset{p_0}{\ket{1}}  \right )  \ket{k^0}\ket{k^1}  U_{c_l}\ket{\hat \tau_{k^0k^1}}  \, , 
\end{equation}
where~\mbox{$\ket{\hat \tau_{k^0k^1}} \equiv \tau_{0, k^0k^1} \ket 0  + \tau_{1, k^0k^1} \ket 1 $}. See Figure~\ref{fig:lcu_basic_b} for the quantum circuit representation of~\eqref{eq:U_hat_Gamma}.

\begin{figure}[t]
\centering
	\includegraphics[scale=0.825]{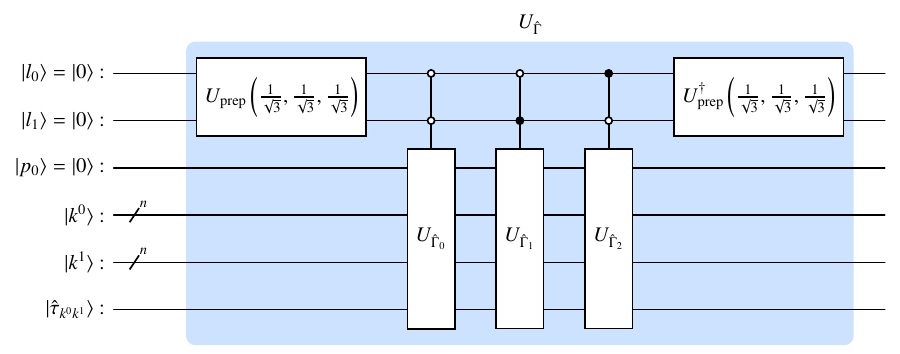} 
\caption{ Quantum circuit for applying the strain Green's operator~$\hat \Gamma_{k^0 k^1}$ using the LCU technique and the decomposition~\eqref{eq:gammaLCU_0}. The discretisation consists of~$2^n \times 2^n$ grid points. In the output, the relevant amplitudes correspond to basis vectors of the form~$\ket 0 \ket 0 \ket 1 \ket{k^0} \ket{k^1}  \ket{\hat \tau_{k^0 k^1}}$, i.e. the ancilla qubits~$\ket{l_0}$ and~$\ket {l_1}$ must be in state~$\ket 0$, and the ancilla qubit~$\ket{p_0}$ in state~$\ket 1$.
\label{fig:lcu_full}}
\end{figure}

The sum of~$U_{\hat \Gamma_l}$ can now be implemented using LCU. As usual, LCU requires a state preparation unitary~$U_{\text{prep}} $ and a select unitary~$U_{\text{select}} $.  According to decomposition~\eqref{eq:gammaLCU_0}, only the unweighted sum of the three unitaries~$U_{\hat \Gamma_l}$ is required. We therefore define the unitaries
\begin{equation} \label{eq:prep_unit}
	U_{\text{prep}} \colon \underset{l_0}{\ket 0} \underset{l_1}{\ket 0}  \mapsto \sum_{l=0}^2 \frac{1}{\sqrt 3} \ket l \,   , 
\end{equation}
and 
\begin{equation}
  	U_{\text{select}} \colon  {\ket l}  \ket{\Phi}  \mapsto  \ket l U_{\hat \Gamma_l} \ket{\Phi} \qquad \text { with }  \quad  \ket \Phi \coloneq  \underset{p_0}{\ket 0} \ket{k^0} \ket {k^1} \ket {\hat \tau_{k^0k^1}} \, . 	
\end{equation}
The linear combination of the unitaries can now be implemented as follows:
\begin{subequations} \label{eq:lcu_gamma}
\begin{align} 
	U_{\text{prep}} \otimes I \colon  &  {\ket 0}^{\otimes 2}  \ket {\Phi}     \mapsto \sum_{l=0}^2 \frac{1}{\sqrt 3} \ket l \ket {\Phi}  \label{eq:lcu_gamma_0} \, ,\\
	U_{\text{select}} \colon  &	\sum_{l=0}^2 \frac{1}{\sqrt 3}  \ket l \ket {\Phi}    \mapsto \sum_{l=0}^2 \frac{1}{\sqrt 3} \ket l U_{\hat  \Gamma_l}\ket  {\Phi}   \label{eq:lcu_gamma_1} \, ,  \\ 
	U^\dagger_{\text{prep}} \otimes I \colon  &  \sum_{l=0}^2 \frac{1}{\sqrt 3} \ket l U_{\hat \Gamma_l}\ket  \Phi   \mapsto  \sum_{l=0}^2 \frac{1}{\sqrt 3} U^\dagger_{\text{prep}}  \ket l U_{\hat \Gamma_l}\ket  \Phi \, .   \label{eq:lcu_gamma_2}
\end{align} 
\end{subequations}
In this section, $I$ is the $2^{2n+2}\times 2^{2n+2}$ identity matrix; elsewhere, its dimension will be clear from the context. In the resulting state, we are only interested in the states when the two qubits~$ \ket {l_0} $ and~$ \ket {l_1} $ representing~$\ket l$ are both in state~$\ket {0}$. These states can be obtained  by applying the operator~\mbox{$(\ket {0} \bra {0} )^{\otimes 2} \otimes I$} to the state in~\eqref{eq:lcu_gamma_2}.  To this end, note that the adjoint~$U_{\text{prep}}^\dagger $ is according to~\eqref{eq:prep_unit} defined as  
\begin{equation}
		U_{\text{prep}}^\dagger \colon {\bra 0}^{\otimes 2} \mapsto \sum_{l=0}^2 \frac{1}{\sqrt 3} \bra l \, .
\end{equation}
Hence, 
\begin{equation}
	( \ket 0 \bra 0 )^{\otimes 2} \otimes I \colon   \sum_{l=0}^2 \frac{1}{\sqrt 3} U^\dagger_{\text{prep}}  \ket l U_{\hat \Gamma_l}\ket  {\Phi} \mapsto 
	  \ket 0^{\otimes 2} \sum_{l=0}^2  \frac{1}{3} U_{\hat \Gamma_l} \ket \Phi   \equiv    \underset{l_0}{\ket 0} \underset{l_1}{\ket 0} \sum_{l=0}^2  \frac{1}{3} U_{\hat \Gamma_l}   \underset{p_0}{\ket 0} \ket{k^0} \ket {k^1} \ket {\hat \tau_{k^0k^1}} \,  . 
\end{equation}
It is straightforward to verify that this is equal to the polarisation stress up to a scaling factor.  To represent the introduced sequence of mappings concisely, we define the composite unitary 
\begin{equation}
U_{\hat \Gamma}  = \left ( U_{\text{prep}}^\dagger \otimes I  \right )  U_{\text{select}} \left ( U_{\text{prep}} \otimes I  \right ) \, .
\end{equation}
The circuit implementation of the derived LCU approach is shown in Figure~\ref{fig:lcu_full}.
%

%
%--------------------------------------------------------------------------------
\subsection{Incremental strain update \label{sec:incrementalUpdate}}
\begin{figure} 
\centering
	\includegraphics[scale=0.78]{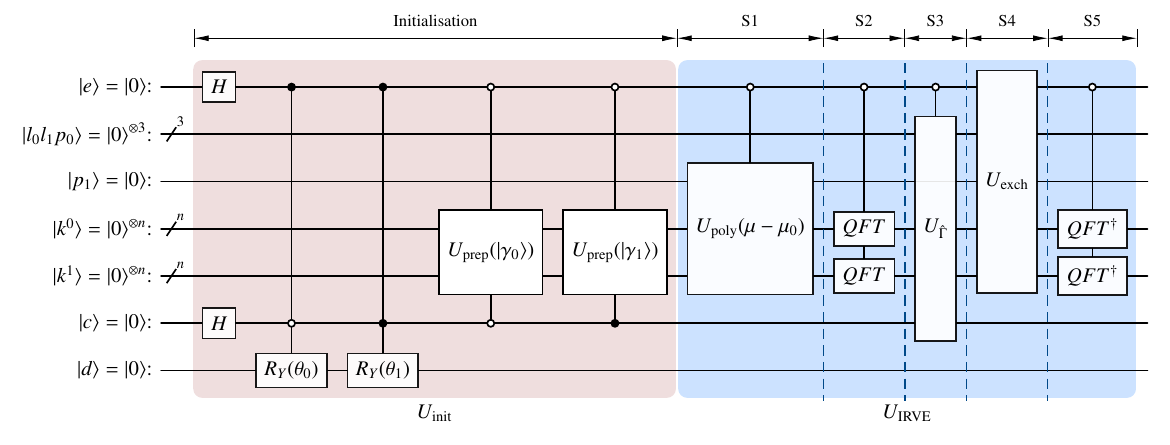} 
\caption{ Quantum circuit for initialising the fixed-point iteration and updating the strain~$\ket{ \gamma_{k^0 k^1}}$.  The gates in the left shaded box serve to initialise the fixed-point iteration.  The two $R_Y$ gates with the arguments~$\theta_0 = 2\cos^{-1} (\bar \gamma_0/N)$ and~$\theta_1 = 2\cos^{-1} (\bar \gamma_1/N)$ encode the prescribed macroscopic strain~$\ket {\overline \gamma}$.  The two other~$U_{\text{prep}}$ gates initialise the components of the initial strain field~$\ket \gamma$.  The gates inside the right shaded box update the strain~$\ket{ \gamma}$. The discretisation consists of~$2^n \times 2^n$ grid points. The relevant amplitudes correspond to basis vectors of the form~\mbox{$\ket 0 \ket {001} \ket 1 \ket{k^0} \ket{k^1}  \ket{\gamma_{k^0k^1}}$}.
\label{fig:circuit_incremental}}
\end{figure}

%--------------------------------------------------------------------------------
%
Next, we introduce the quantum circuit for computing the strain vector~\mbox{${\ket \gamma}^{(s+1)}  \in \mathbb R^{2 N^2} $ at iteration step $(s+1)$} for the given strain vector~\mbox{${\ket \gamma}^{(s)} \in \mathbb R^{2 N^2}$}  at step~$(s)$. As before, we omit the iteration step index~$(s)$ to keep the notation simple. The quantum circuit for the proposed unitary~$U_{\text{IRVE}}$  is shown in~Figure~\ref{fig:circuit_incremental}. The shown circuit also includes the encoding of the macroscopic strain vector~$\ket {\overline \gamma}$  and the encoding of an initial strain field vector~$\ket{\gamma}$ at the beginning of the iteration. 

The full state vector of the circuit is formed by the following qubits and registers:  
\begin{equation} 
	\ket {e} \ket {l_0 l_1 p_0} \ket{p_1} \ket {k^0}  \ket{k^1} \ket c \ket d  \,  , 
\end{equation}
where~$e, l_0, l_1, p_0, p_1, c, d \in \{0, \,1 \} $ and $k^0, k^1 \in \{0, \, 1, \, \dotsc, \, N-1 \}$. The registers~$\ket{k^0}$ and~$\ket{k^1}$ consist each of~$n = \log_2 N$ qubits and represent the indices of the grid points. As in previous Section~\ref{sec:incrementalFourierSpace}, the ancilla qubits~$\ket{l_0 l_1 p_0}$ belong to the Fourier-space solution unitary~$U_{\hat \Gamma}$ mapping the stresses to strains. The new ancilla qubit~$\ket {p_1}$ is used while mapping the strains to stresses. The remaining qubits~$\ket e$, $\ket c$, and~$\ket d$ are for encoding the strain and stress fields, and the prescribed macroscopic strain. 

To begin with, all qubits are in the state~$\ket 0$ and must be initialised prior to applying the~$U_{\text{IRVE}}$. Without loss of generality, we assume that the two components of the initial strain vector~$\ket{\gamma}$ over all grid points are independently normalised, whereas the prescribed macroscopic strain~$\ket {\overline{\gamma}}$ need not be normalised. According to Figure~\ref{fig:circuit_incremental}, the initial quantum state vector is constructed as follows. The application of the first two Hadamard gates yields the state
\begin{equation}
	 \underset{e}{\ket 0}   \underset{l_0 l_1 p_0 }{\ket {000}}   \underset{p_1}{\ket{0}}  \underset{k^0}{\ket {0}}  \underset{k^1}{\ket {0}}  \underset{c}{\ket {0}}  \underset{d}{\ket {0}}  \mapsto   \frac{1}{2}  \left  ( \underset{e}{\ket 0} + \underset{e}{\ket 1} \right ) \underset{l_0 l_1 p_0}{\ket {000}}  \underset{p_1}{\ket{ 0}}  \underset{k^0}{\ket {0}}  \underset{k^1}{\ket {0}}  \left (\underset{c}{\ket 0} + \underset{c}{\ket {1}} \right )    \underset{d}{\ket {0}}\, .
\end{equation}
The two components of the prescribed macroscopic strain~$\ket {\overline \gamma}$ are encoded using the two controlled rotation gates 
\begin{equation}
	R_Y(2 \cos^{-1} (\overline \gamma_0/N)) = \frac{1}{N} \begin{pmatrix}  \overline \gamma_0 & - \sqrt{N^2- \overline \gamma_0} \\ \sqrt{N^2- \overline \gamma_0}  & \overline \gamma_0 \end{pmatrix} \,  , \quad  R_Y(2 \cos^{-1} (\overline \gamma_1/N)) = \frac{1}{N} \begin{pmatrix}  \overline \gamma_1 & - \sqrt{N^2- \overline \gamma_1} \\ \sqrt{N^2- \overline \gamma_1}  & \overline \gamma_1 \end{pmatrix}  \, .
\end{equation}
The need for the factor~$1/N$ will become evident in the following. Similarly, the two components of the strain field vector~$\ket \gamma$ are encoded using the controlled state preparation unitaries~$U_{\text {prep}} (\ket {\gamma_0})$ and~$U_{\text{prep}}(\ket{ \gamma_1})$. The ancilla qubit~\mbox{$\ket c \in \{ \ket 0, \, \ket 1 \}$} is used to distinguish between the two strain components, and the ancilla qubit~\mbox{$\ket e \in \{ \ket 0, \, \ket 1 \}$} to distinguish between the strain field and the prescribed macroscopic FE strain. The resulting state vector reads
\begin{equation}
\begin{aligned}
	 & \frac{1}{2} \sum_{k^0, \, k^1=0}^{N-1} \underset{e} {\ket 0}  \ket{000}  \ket 0  \ket {k^0} \ket {k^1} \left ( \gamma_{0, \, k^0k^1} \underset{c} {\ket 0} + \gamma_{1, \, k^0k^1} \underset{c} { \ket 1} \right ) \ket 0   \\
	 + &\frac{1}{2 N} \underset{e} {\ket 1}  \ket{000}  \ket 0  \underset{k^0} {\ket 0} \underset{k^1}{\ket 0} \underset{c}{\ket 0} \left ( \overline {\gamma}_0  \underset{d} {\ket 0} + \sqrt{N^2 -\overline {\gamma}_0^2}  \underset{d} {\ket 1} \right ) + \frac{1}{2 N} \underset{e}{\ket 1}  \ket{000}  \ket 0  \underset{ k^0} {\ket 0} \underset{k^1}{\ket 0}  \underset{c}{\ket 1} \left ( \overline {\gamma}_1  \underset{d}{\ket 0} + \sqrt{N^2 -\overline {\gamma}_1^2}  \underset{d}{\ket 1} \right ) \, .
\end{aligned}
\end{equation}
The states with~$\ket{d} = \ket{1}$ have only been introduced to preserve normalisation and are not required in the following. Henceforth, we always discard the unnecessary states and abbreviate the state vector as 
\begin{equation}
\begin{aligned}
	 & \frac{1}{2} \sum_{k^0, \, k^1 = 0}^{N-1} \underset{e}{\ket 0}  \ket{000}  \ket 0  \ket {k^0} \ket {k^1} \left (\gamma_{0, \, k^0k^1} \underset{c}{\ket 0} + \gamma_{1, \, k^0k^1} \underset{c} {\ket 1} \right ) \underset{d}{\ket 0}   \\
	 + & \frac{\overline {\gamma}_0}{2 N} \underset{e}{\ket 1}  \ket{000}  \ket 0 \underset{k^0} {\ket 0} \underset{k^1}{\ket 0}   \underset{c}{\ket 0} \underset{d}{\ket 0}  +  \frac{\overline {\gamma}_1}{2 N}  \underset{e}{\ket 1}  \ket{000}  \ket 0  \underset{k^0} {\ket 0} \underset{k^1}{\ket 0}   \underset{c}{\ket 1} \underset{d}{\ket 0}  + \dotsc \, .
\end{aligned}
\end{equation}

With the preparation of the initial state vector at iteration step~$s=0$ completed, we proceed to the discussion of the unitary~$U_{\text{IRVE}}$, cf. Figure~\ref{fig:circuit_incremental}. The sequence of operations in the circuit corresponds to the five steps summarised in Box~\ref{tb:rve_problem}. In S1,  the polarisation stress field~$\ket \tau \in \mathbb R^{2 N^2}$ is computed by applying the controlled unitary~$U_{\text{poly}} (\mu (k^0, k^1)-\mu_0)$  yielding the state vector
\begin{equation}
\begin{aligned}
	 & \frac{1}{2} \sum_{k^0, \, k^1=0}^{N-1}  \underset{e}{\ket 0}  \ket{000}  \underset{p_1}{\ket 1}  \ket {k^0} \ket {k^1} \left (\tau_{0, \, k^0k^1} \underset{c}{\ket 0} + \tau_{1, \, k^0k^1} \underset{c} {\ket 1} \right ) \underset{d}{\ket 0}   \\
	 + & \frac{\overline {\gamma}_0}{2 N} \underset{e}{\ket 1}  \ket{000}  \underset{p_1}{\ket 0} \underset{k^0} {\ket 0} \underset{k^1}{\ket 0}   \underset{c}{\ket 0} \underset{d}{\ket 0}  +  \frac{\overline {\gamma}_1}{2 N}  \underset{e}{\ket 1}  \ket{000}  \underset{p_1}{\ket 0}\underset{k^0} {\ket 0} \underset{k^1}{\ket 0}   \underset{c}{\ket 1} \underset{d}{\ket 0}  + \dotsc \, .
\end{aligned}
\end{equation}
The componentwise multiplication of the strain field~$\ket{\gamma}$ by the shear modulus field~\mbox{$(\mu (k^0, k^1)-\mu_0)$} is non-length-preserving and requires the ancilla~$\ket {p_1}$. For~$s= 0$, the states with~$\ket {p_1} = \ket{1}$ represent the physical component of the polarisation stress and the omitted states with~$\ket{p_1} = \ket{0}$ the complementary component. 
Subsequently, in S2, the application of the two controlled one-dimensional QFTs yields the partially Fourier-transformed state  
\begin{equation}
\begin{aligned}
	 & \frac{1}{2 N} \sum_{k^0, \, k^1=0}^{N-1}  \underset{e}{\ket 0}  \ket{000}  \underset{p_1}{\ket 1}  \ket {k^0} \ket {k^1} \left (\hat \tau_{0, \, k^0k^1} \underset{c}{\ket 0} + \hat \tau_{1, \, k^0k^1} \underset{c} {\ket 1} \right ) \underset{d}{\ket 0}   \\
	 + & \frac{\overline {\gamma}_0}{2 N} \underset{e}{\ket 1}  \ket{000}  \underset{p_1}{\ket 0} \underset{k^0} {\ket 0} \underset{k^1}{\ket 0}   \underset{c}{\ket 0} \underset{d}{\ket 0}  +  \frac{\overline {\gamma}_1}{2 N}  \underset{e}{\ket 1}  \ket{000}  \underset{p_1}{\ket 0} \underset{k^0} {\ket 0} \underset{k^1}{\ket 0}  \underset{c}{\ket 1} \underset{d}{\ket 0}  + \dotsc \, .
\end{aligned}
\end{equation}
The updated strain field~$\ket{\hat  \gamma}$ in S3 is computed using the controlled unitary~$U_{\hat \Gamma}$. As introduced in Section~\ref{sec:incrementalFourierSpace},~$U_{\hat \Gamma}$ depends on the ancilla qubits~$\ket{l_0l_1p_0}$ and the physical components in the output state correspond to state~\mbox{$\ket{l_0l_1p_0} = \ket{001}$}. Hence, after the application of~$U_{\hat \Gamma}$ the state vector reads
\begin{equation}
\begin{aligned}
	 & \frac{1}{2 N} \sum_{k^0, \, k^1=0}^{N-1}  \underset{e}{\ket 0} \underset{l_0 l_1 p_0} {\ket{001}}  \underset{p_1}{\ket 1}  \ket {k^0} \ket {k^1} \left (\hat \gamma_{0, \, k^0k^1} \underset{c}{\ket 0} + \hat \gamma_{1, \, k^0k^1} \underset{c} {\ket 1} \right ) \underset{d}{\ket 0}   \\
	 + & \frac{\overline {\gamma}_0}{2 N} \underset{e}{\ket 1}   \underset{l_0 l_1 p_0}{ \ket{000}} \underset{p_1}{\ket 0} \underset{k^0} {\ket 0} \underset{k^1}{\ket 0}   \underset{c}{\ket 0} \underset{d}{\ket 0}  +  \frac{\overline {\gamma}_1}{2 N}  \underset{e}{\ket 1}  \underset{l_0 l_1 p_0}{\ket{000}}  \underset{p_1}{\ket 0}  \underset{k^0} {\ket 0} \underset{k^1}{\ket 0}   \underset{c}{\ket 1} \underset{d}{\ket 0}  + \dotsc \, .
\end{aligned}
\end{equation}
In S4, the macroscopic strain is applied by setting~$\hat \gamma_{0,00}= \overline \gamma_0$ and~$\hat \gamma_{1,00}= \overline \gamma_1$  for~$k^0=k^1=0$.  We implement this using the unitary~$U_{\text{exch}}$, which exchanges, i.e. permutes, the relevant components of the state vector so that 
\begin{equation}
\begin{aligned} \label{eq:exch_bar_gamma}
	 & \frac{1}{2 N}  \underset{e}{\ket 0} \underset{l_0 l_1 p_0} {\ket{001}}  \underset{p_1}{\ket 1}  \underset{k^0}{\ket {0}} \underset{k^1}{\ket {0}} \left (\overline \gamma_{0} \underset{c}{\ket 0} + \overline \gamma_{1} \underset{c} {\ket 1} \right ) \underset{d}{\ket 0}  + \frac{1}{2 N} \sum_{\substack{k^0, \, k^1=0 \\ (k^0, \, k^1)\neq(0,0)}}^{N-1} \underset{e}{\ket 0} \underset{l_0 l_1 p_0} {\ket{001}}  \underset{p_1}{\ket 1}  \ket {k^0} \ket {k^1} \left (\hat \gamma_{0, \, k^0k^1} \underset{c}{\ket 0} + \hat \gamma_{1, \, k^0k^1} \underset{c} {\ket 1} \right ) \underset{d}{\ket 0}   \\
	 + & \frac{\hat {\gamma}_{0,00}}{2 N} \underset{e}{\ket 1}   \underset{l_0 l_1 p_0}{ \ket{000}}  \underset{p_1}{\ket 0}  \underset{k^0}{\ket 0} \underset{k^1}{\ket 0}  \underset{c}{\ket 0} \underset{d}{\ket 0}  +  \frac{\hat {\gamma}_{1,00}}{2N}  \underset{e}{\ket 1}  \underset{l_0 l_1 p_0}{\ket{000}}  \underset{p_1}{\ket 0}  \underset{k^0}{\ket 0} \underset{k^1}{\ket 0}  \underset{c}{\ket 1} \underset{d}{\ket 0}  + \dotsc \, .
\end{aligned}
\end{equation}
See~\ref{app:graycode_exch}  for an efficient implementation of~$U_{\text{exch}}$. Finally, in S5, we apply the inverse QFT to obtain the updated strain field with the state vector 
\begin{equation}
\begin{aligned}
 \frac{1}{2} \sum_{k^0, \, k^1 = 0}^{N-1}  \underset{e}{\ket 0} \underset{l_0 l_1 p_0} {\ket{001}}  \underset{p_1}{\ket 1}  \ket {k^0} \ket {k^1} \left (\gamma_{0, \, k^0k^1} \underset{c}{\ket 0} +  \gamma_{1, \, k^0k^1} \underset{c} {\ket 1} \right ) \underset{d}{\ket 0}  + \dotsc \, .
\end{aligned}
\end{equation}
In the entire state vector, only the components corresponding to ancilla states~$\ket e = \ket 0$,  $\ket {l_0 l_1 p_0} = \ket{001}$, $\ket {p_1} = \ket 1$ and  $\ket d = \ket 0$ represent physical strains; the other components are required to implement the non-unitary operations. It is worth emphasising that the total number of ancilla qubits per iteration step is independent of the grid size~$N^2$.

%
%--------------------------------------------------------------------------------
\subsection{Fixed-point iteration \label{sec:fixedPointIter}}
\begin{figure} 
	\centering
	\includegraphics[scale=0.8]{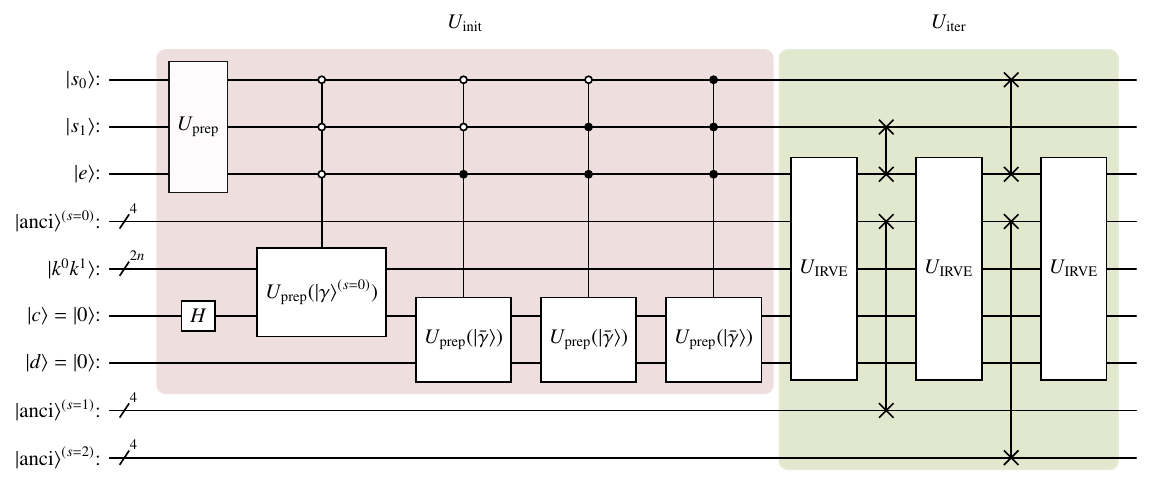} 
	\caption{Quantum circuit for fixed-point iteration with~$S=3$ iteration steps. \label{fig:fixed_point_iteration}}
\end{figure}

%--------------------------------------------------------------------------------
%
We are now ready to introduce the quantum algorithm for the fixed-point iteration comprising several iteration steps (Figure~\ref{fig:fixed_point_iteration}). We assume that the total number of iteration steps~$S$ is fixed from the outset. The proposed circuit entails~$S$ back-to-back applications of the~$U_{\text{IRVE}}$ circuit introduced in the previous section. The initialisation step introduced earlier is extended because copying or deleting components of the state vector is impossible (the no-cloning theorem). Specifically, during the initialisation step, we encode~$S$ copies of the macroscopic FE strain~$\ket{\overline \gamma}$ to enforce the macroscopic strain within every iteration step~$s \in \{0, \, \dotsc , \,  S-1\}$ by setting the zeroth Fourier components of the strain vector~$\hat\gamma_{0,00}^{(s)} = \overline \gamma_0$ and~$\hat \gamma_{1,00}^{(s)} = \overline \gamma_1$. Furthermore, the ancillary qubits required for LCU and polynomial interpolation, i.e.~$\ket{l_0 l_1 p_0 p_1}$ denoted as~\mbox{$\ket{\text {anci}} =  \ket{l_0 l_1 p_0 p_1}$} in Figure~\ref{fig:fixed_point_iteration} cannot be reused. Consequently, a new set of initialised ancilla qubits in state~$\ket{0000}$,  must be provided for each iteration step.

In the proposed quantum circuit in Figure~\ref{fig:fixed_point_iteration}, the register~$\ket {s_0, \,  \dotsc, \,  s_{S-2}}$, where~$S=3$, in conjunction with the previously introduced~$\ket e$ qubit,  is used to keep track of the number of times~$U_{\text{IRVE}}$ is applied.  The generalisation to arbitrary~$S$ is straightforward. Using a standard state preparation unitary~$U_{\text{prep}}$ the state vector corresponding  to qubits~$\ket {s_0  s_1} \ket e$ is initialised as follows 
\begin{equation} \label{eq:i0i1s_state_fixed_point}
	\ket {00} \ket{0} \mapsto \frac{1}{2} \left ( \ket {00} \ket{0} + \ket {00} \ket{1} + \ket {01} \ket{1} + \ket {11} \ket{1} \right ) \, .
\end{equation}
These components facilitate the encoding of the initial strain~$\ket{\gamma}^{(s=0)}$ and the three copies of the prescribed strain~$\ket{\overline \gamma}$ into the four implied subspaces. Each subspace is selected by suitably choosing the control states for the respective state preparation unitaries~$U_{\text{prep}}$. Although we denote all unitaries for encoding a given vector as the amplitudes of a quantum state as~$U_{\text{prep}}$,  their specific implementation depends on the vector to be encoded. For instance, as discussed in the previous section, each of the unitaries for encoding~$\ket{\overline \gamma}$ consists of two~$R_Y$ gates, while the unitary for encoding~$\ket{\gamma}$ is composed of two unitaries for encoding its components.  

According to Figure~\ref{fig:fixed_point_iteration}, for all three copies of the prescribed strain~$\ket{\overline \gamma}$ the~$\ket e$ qubit is in state~$\ket 1$ and for the initial strain field~$\ket{\gamma}$ it is in state~$\ket 0$. The qubits~$\ket{s_0 s_1}$ are in states~$\ket {00}$, $\ket {01}$ and~$\ket{11}$ for the three copies of the prescribed strain.  As discussed in the previous section, the unitary~$U_{\text{IRVE}}$ assumes that for~$\ket{\overline \gamma}$ the qubit~$\ket e$ is in state~$\ket 1$ and for~$\ket{\gamma}$ it is in state~$\ket 0$.  Note that after each application of~$U_{\text{IRVE}}$, the initial strains~$\ket{\overline \gamma}$ are assigned to states with~$ \ket{e} = \ket{0}$ due to quantum entanglement. The specific choice of the states in~\eqref{eq:i0i1s_state_fixed_point}  and the swap gates after each application of~$U_{\text{IRVE}}$  in Figure~\ref{fig:fixed_point_iteration} ensure that the relevant initial strains~$\ket{\overline \gamma}$ are reassigned to states with~$\ket{e} = \ket{1}$. Furthermore, as indicated in Figure~\ref{fig:fixed_point_iteration}, after every application of~$U_{\text{IRVE}}$, the set of used ancilla qubits~$\ket{\text{anci}}^{(s=0)}$ are replaced with the initialised ancilla qubits in state~$\ket{\text{anci}}^{(s)} = \ket{0}^{\otimes 4}$ using the swap gates. At the end of the iteration after $S$ applications of~$U_{\text{IRVE}}$  the relevant strain vector~$\ket{\gamma}$ corresponds to states with~$\ket{e} = \ket{ 0}$.

%
%--------------------------------------------------------------------------------
\section{Quantum parallelisation \label{sec:quantum_parallel}}
%--------------------------------------------------------------------------------
%
Quantum computing offers an exponential speedup in evaluating a single RVE. However, the advantage of quantum computing becomes even more pronounced when evaluating an ensemble of RVEs corresponding to all quadrature points of a large FE model. Owing to superposition and entanglement, all the RVEs in the ensemble can be processed in parallel. We begin by illustrating the basic idea of quantum parallelisation in Section~\ref{sec:quantum_parallel_motivating} using the quantum Fourier transform in a linear-algebraic setting. The parallelisation of RVE problems is then introduced in Section~\ref{sec:quantum_parallel_rve_update}.

%
%--------------------------------------------------------------------------------
\subsection{Motivating example \label{sec:quantum_parallel_motivating}}
\begin{figure}
\centering
\includegraphics[scale=0.875]{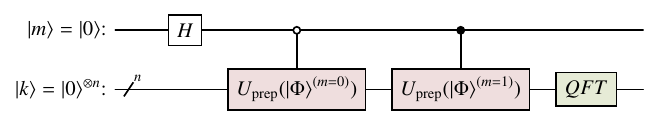} 
\caption{Quantum Fourier transformation of two vectors $\ket{\Phi}^{(m=0)}$ and $\ket{\Phi}^{(m=1)}$ utilising a single~$QFT$ unitary. \label{circ:parallel_qft}}
\end{figure}

%--------------------------------------------------------------------------------
%
We consider the quantum Fourier transformation of~$M$ vectors~$ \ket {\Phi}^{(m=0)},  \ket {\Phi}^{(m=1)}, \dotsc ,  \ket {\Phi}^{(m=M-1)} \in \mathbb C^N$. The number of qubits required to encode~$M$ and~$N$ are~\mbox{$m = \log_2 M$} and~\mbox{$n = \log_2 N$}, respectively. We assume that~$M$ and~$N$ are powers of~$2$.  As shown in Figure~\ref{circ:parallel_qft} for the case $M=2$, the Fourier transforms $\ket {\hat \Phi}^{(m=0)},  \ket {\hat \Phi}^{(m=1)} \in \mathbb C^N$ can be computed in parallel using a single unitary $QFT \in \mathbb C^{N\times N}$.  In the depicted circuit, the two vectors are encoded using the state preparation unitaries~\mbox{$U_{\text{prep}} (\ket{\Phi}^{(m=0)}) , U_{\text{prep}} (\ket{\Phi}^{(m=1)}) \in \mathbb C^{N\times N}$}, although the specific method of state preparation is immaterial for the present discussion. 

It is instructive to examine the sequence of unitary transformations implemented by the circuit. The state vector is formed by the registers~$\ket m$ and~$\ket k$, where~$m \in \{0,\, 1\}$ and~$k \in \{0, \, 1, \, \dotsc, \, N-1 \}$. The application of the Hadamard gate~$H$ on the~$\ket m$ qubit and the subsequent application of the controlled state preparation unitaries $U_{\text{prep}} (\ket{ \Phi}^{(m=0)} )$ and $U_{\text{prep}} (\ket{ \Phi}^{(m=1)} ) $ yield 
\begin{equation}
	\ket 0 \ket 0^{\otimes n} \mapsto \frac{1}{\sqrt{2}} \ket 0   \ket { \Phi}^{(m=0)}   + \frac{1}{\sqrt{2}}  \ket 1 \ket { \Phi}^{(m=1)} \, .
\end{equation}
A single~$QFT \in \mathbb C^{N\times N}$ is applied to the~$\ket k$ register, acting on both components of the superposition to produce the final state
\begin{equation} \label{eq:mot_ex_final_state}
	 \frac{1}{\sqrt{2}} \ket 0   \ket {\hat  \Phi}^{(m=0)}   + \frac{1}{\sqrt{2}}  \ket 1 \ket {\hat \Phi}^{(m=1)} \, .
\end{equation}
It is also useful to consider the unitary implemented by the entire circuit. In block-matrix form, this unitary can be written as the product 
\begin{equation}
\frac{1}{\sqrt{2}} 
\begin{pmatrix*}[c]
	QFT & \\ 
	& QFT
\end{pmatrix*}
\begin{pmatrix*}[c]
	I & \\ 
	&U_{\text{prep}}\left (\ket{ \Phi}^{(m=1)} \right )  
\end{pmatrix*}
\begin{pmatrix*}[c]
	U_{\text{prep}} \left (\ket{ \Phi}^{(m=0)} \right )  & \\ 
	& I
\end{pmatrix*}
\begin{pmatrix*}[r]
	I  &  I \\ 
	I & -I  
\end{pmatrix*} \, ,
\end{equation}
where~$I$ is the identity matrix of dimension~$2^n \times 2^n$. Acting on the initial state~$\ket 0 \ket 0^{\otimes n} = \begin{pmatrix}  1 \cdot \ket 0^{\otimes n} &  0 \cdot \ket 0^{\otimes n} \end{pmatrix}^\trans$, this unitary  yields the Fourier-transformed state~\eqref{eq:mot_ex_final_state}.

\begin{figure}
\centering
\includegraphics[scale=0.9]{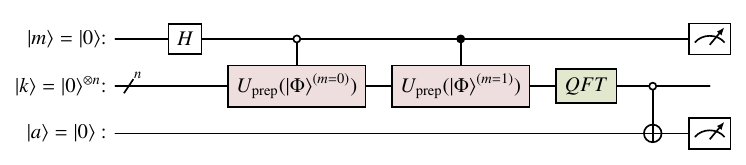} 
\caption{Quantum Fourier transformation of two vectors $\ket{\Phi}^{(m=0)}$ and $\ket{\Phi}^{(m=1)}$ utilising a single~$QFT$ unitary, followed by the measurement of the zeroth Fourier components of $\ket{\hat{\Phi}}^{(m=0)}$ and $\ket{\hat{\Phi}}^{(m=1)}$. \label{circ:parallel_qft_mean}}
\end{figure}

In homogenisation problems, we are interested only in the average stress of the RVE. The average of a vector can be determined, for instance, by computing its Fourier transform and measuring its zeroth component.  As an example, the zeroth components of the vectors~$\ket{\hat{\Phi}}^{(m=0)}$ and~$\ket{\hat{\Phi}}^{(m=1)}$ can be measured using the circuit shown in Figure~\ref{circ:parallel_qft_mean}. Compared to the original circuit in Figure~\ref{circ:parallel_qft}, it contains one ancilla qubit~$\ket a$ and a multi-controlled $\mathit{NOT}$ gate.   
Prior to the measurement of the~$\ket m$ and~$\ket a$ qubits, the state vector of the circuit is given by 
\begin{equation}
	\ket \Psi = \frac{1}{\sqrt{2}} \hat  \Phi_0^{(m=0)}  \ket 0 \ket 0^{\otimes n} \ket 1 + \frac{1}{\sqrt{2}}  \hat \Phi_0^{(m=1)}  \ket 1 \ket 0^{\otimes n} \ket 1  + \frac{1}{\sqrt{2}} \sum_{k=1}^{N-1}  \left ( \hat \Phi_k^{(m=0)} \ket 0   +   \hat \Phi_k^{(m=1)} \ket 1 \right ) \ket k \ket 0 \, .
\end{equation}
To determine~$ \hat  \Phi_0^{(m=0)} $ and~$ \hat  \Phi_0^{(m=1)} $ via measurement, we define the four projectors
\begin{equation}
	\Pi_0 = \ket 0 \bra 0 \otimes I \otimes  \ket 1 \bra 1 \,  , \quad \Pi_1 = \ket 1 \bra 1 \otimes I \otimes \ket 1 \bra 1  \,  , \quad \Pi_2 = \ket 0 \bra 0 \otimes  I \otimes \ket 0 \bra 0 \,    , \quad \Pi_3 = \ket 1 \bra 1  \otimes I  \otimes \ket 0 \bra 0 \, .
\end{equation}
As required,~$\sum_{j} \Pi_j =  I$, where $I$ is the \mbox{$2^{n+2}\times 2^{n+2}$} identity matrix, and~$\Pi_j  \Pi_k = 0$ if~$j \neq k$. These projection matrices are applied by measuring the~$\ket m$ and~$\ket a$ qubits, causing the quantum state~$\ket \Psi$ to collapse to one of the four states 
\begin{equation}
 \frac{ \Pi_j \ket \Psi}{ \sqrt{\bra \Psi  \Pi_j \ket \Psi}}  \, .
\end{equation}
The probability of finding the system in state~$j$ is given by 
\begin{equation}
	p(j) = \bra \Psi  \Pi_j \ket \Psi \, . 
\end{equation}
Consequently, $ \hat  \Phi_0^{(m=0)} $ and~$ \hat  \Phi_0^{(m=1)} $ correspond to the probabilities~$p(0)$ and~$p(1)$, respectively.  

%
%
%--------------------------------------------------------------------------------
\subsection{Simultaneous solution of all RVEs \label{sec:quantum_parallel_rve_update}}
%--------------------------------------------------------------------------------
%
To explain the basic construction of the \emph{QAFE$^2$} framework, we assume that all RVE problems, including their material properties and discretisation, are identical.  The RVE problems differ only in terms of the prescribed macroscopic FE strain~$\ket{\overline {\gamma}}$. The number of RVE problems is~$M$ and their prescribed macroscopic FE strains are~$\ket{\overline{\gamma}}^{(m)}$, where~$m \in \{0, \, 1, \, \dotsc, \, M-1 \}$.

\begin{figure} 
	\centering
	\includegraphics[scale=0.9]{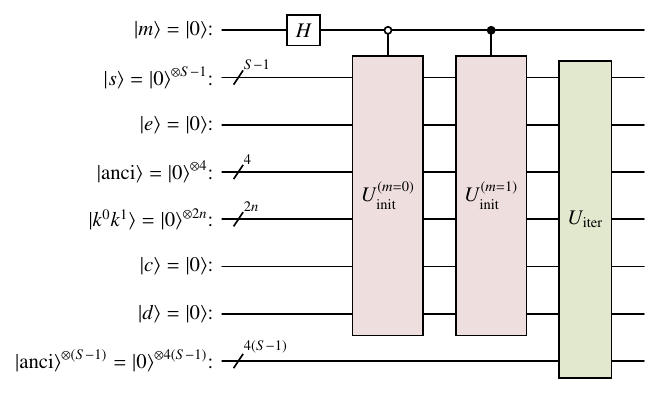} 
	\caption{Quantum circuit for the simultaneous solution of two ($M=2$) RVE problems. The controlled unitaries~$U_{\text{init}}^{(m=0)}$ and~$U_{\text{init}}^{(m=1)}$ encode the prescribed initial strain~$\ket{\gamma}^{(s=0)}$ and the macroscopic FE strains~$\ket{\overline{ \gamma}}^{(m=0)}$ and~$\ket{\overline{\gamma}}^{(m=1)}$, and the unitary~$U_{\text{iter}}$ implements the fixed-point iteration.  \label{fig:qafe2_2}}
\end{figure}

The implementation of \emph{QAFE$^2$} follows the motivating example in the previous section. The quantum circuit for \emph{QAFE$^2$}  for the case~$M=2$, involving two RVE problems, is depicted in Figure~\ref{fig:qafe2_2}. Its strong similarity to the circuit in Figure~\ref{circ:parallel_qft} of the linear-algebraic motivating example is evident. The controlled unitaries~$U_{\text{init}}^{(m=0)}$ and~$U_{\text{init}}^{(m=1)}$ for initialisation and the unitary~$U_{\text{iter}}$ for the fixed-point iteration, are exactly the same as those used for single RVEs introduced in Section~\ref{sec:fixedPointIter}. Crucially, there are two unitaries for initialisation that encode the prescribed initial strain~$\ket{\gamma}^{(s=0)}$ and the macroscopic FE strains~$\ket{\overline{ \gamma}}^{(m=0)}$ and~$\ket{\overline{ \gamma}}^{(m=1)}$, but only a single unitary for the fixed-point iteration. The Hadamard gate~$H$ acting on the~$\ket {m}$ register creates an equal superposition state~\mbox{$\ket m = 1/\sqrt{2} (\ket 0 + \ket 1)$}, and the unitaries for initialisation are conditioned either on~$\ket m = \ket 0$ or~$\ket m = \ket 1$. Intuitively, the quantum state is partitioned into two orthogonal subspaces: the first RVE problem with~$m=0$ is assigned to the subspace with~$\ket m = \ket 0$, and the second RVE problem with~$m=1$ is assigned to the subspace with~$\ket m = \ket 1$.  Owing to quantum entanglement, a single instance of the fixed-point iteration unitary~$U_{\text{iter}}$ is sufficient to implement the fixed-point iteration in both subspaces. 

\begin{figure} 
	\centering
	\includegraphics[scale=0.9]{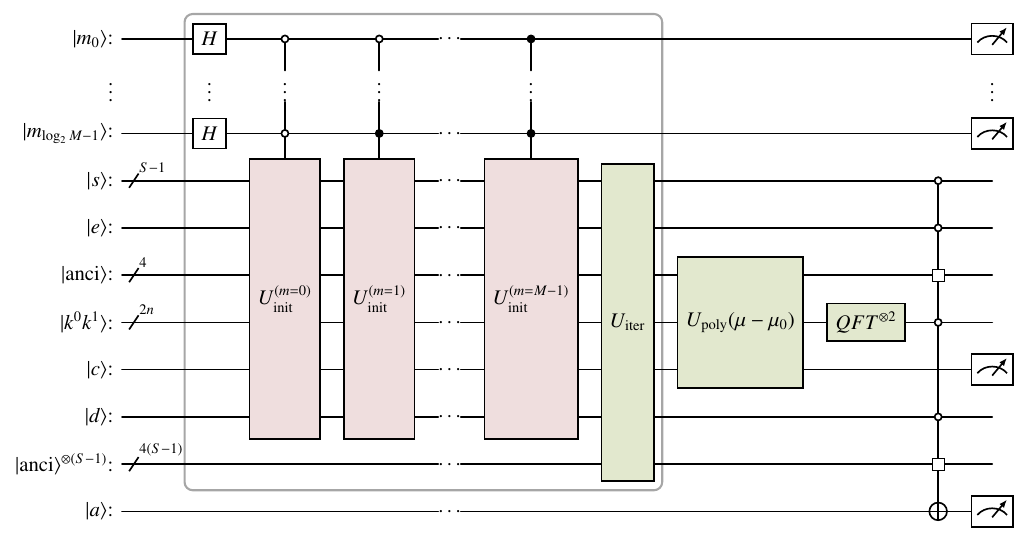} 
	\caption{ Quantum circuit for the simultaneous solution (large box) and measurement of average stress in the~$M$ RVEs.  The multi-controlled unitaries~$U_{\text{init}}^{(m)}$  encode the prescribed initial strain~$\ket{\gamma}^{(s=0)}$ and the macroscopic FE strains~$\ket{\overline{\gamma}}^{(m)}$, where~$m \in \{0, \, 1, \, \dotsc, \, M-1 \}$.  A single fixed-point iteration unitary~$U_{\text{iter}}$ is sufficient to update the strains in all RVEs simultaneously. Prior to measuring the average stresses~$\ket{\overline{{\sigma}}}^{(m)}$, the Fourier transformed stresses are determined by applying the unitary $U_{\text{poly}}$ together with the unitary $QFT ^{\otimes 2}$. \label{fig:qafe2_M}}
\end{figure}

The extension of \emph{QAFE$^2$} to the case with~$M > 2$ is shown in Figure~\ref{fig:qafe2_M}. Assuming that~$M$ is a power of~$2$, an enlarged register~$\ket{m_0  m_{1}  \dotsc m_{\log_2 M-1} } $ is used to create the equal superposition state. This implies a binary enumeration of the RVE problems, which are indexed by~\mbox{$ m_0, m_{1}, \dotsc, m_{\log_2 M-1}  \in \{ 0, 1\}$}. The initial set of Hadamard gates~$H^{\otimes \log_2 M}$ in the circuit create the equal superposition state 
\begin{equation}
	 \frac{1}{\sqrt{M}} \sum_{m_0=0}^{1}  \sum_{m_{1}=0}^{1}   \cdots  \sum_{m_{\log_2 M-1}=0}^{1} \ket{m_0 m_{1}  \dotsc  m_{\log_2 M-1} }  \, .
\end{equation}
The subsequent multi-controlled unitaries $U_{\text{init}}^{(m)}$ are again the same as the one introduced in Section~\ref{sec:fixedPointIter}. They differ only in terms of the encoded prescribed macroscopic strain~$\ket{\overline{ \gamma}}^{(m)}$. The indicated controls of~$U_{\text{init}}^{(m)}$ select the relevant orthogonal subspace.  After this initialisation step, a single instance of the~$U_{\text{iter}}$ unitary suffices to apply the fixed-point iteration simultaneously across all~$M$ subspaces. The average stresses~$\ket{\overline{{\sigma}}}^{(m)}$ are determined by simultaneously computing the Fourier transform of the stress field in all subspaces. The final multi-controlled $\mathit{CNOT}$ gate selects the appropriate subspaces associated with the LCU construction, the polynomial approximation and the fixed-point iteration. The empty control circles indicate conditioning on all qubits in the corresponding register being in the state~$\ket{0}$. The box-shaped controls take into account our choice of state~$\ket{1}$ as the relevant subspace for the polynomial interpolation.

Although we have not yet focused on the computational complexity of the introduced algorithms, it is worth emphasising that the overhead in \emph{QAFE$^2$} arises exclusively from encoding the initial state. There is only a single unitary to implement the fixed-point iteration. As will be numerically demonstrated, the complexity of the \emph{QAFE$^2$} depends approximately linearly on~$M$ with a mild overhead due to the implementation of the multi-controlled gates required to select an appropriate subspace.

%
%--------------------------------------------------------------------------------
\section{Examples \label{sec:examples}}
%--------------------------------------------------------------------------------
%
We present a series of numerical experiments to assess the computational complexity and accuracy of the proposed \emph{QAFE$^2$} framework. All circuits are implemented in Qiskit and executed on a noiseless state vector simulator~\cite{qiskit2024}. We assess computational complexity by expressing the circuits in terms of only two-qubit $\mathit{CNOT}$ and single-qubit~$U_3$ rotation gates. See~\cite{liu2024towards,febrianto2025quantum} for the definition of~$U_3$. The gate set $\{ \mathit{CNOT}, U_3\}$ is universal, meaning that any quantum circuit can be expressed using these two gates. We first study the solution of a one-dimensional RVE problem with a known analytical solution and with prescribed single and multiple macroscopic FE strains. Subsequently, we consider a two-dimensional RVE problem with a known analytical solution and a single macroscopic FE strain.  
%
%--------------------------------------------------------------------------------
\subsection{One-dimensional RVE \label{sec:example_1d}}
%--------------------------------------------------------------------------------
%
We consider a one-dimensional RVE defined on the domain $\Omega=(0, \, L)$ with periodic boundary conditions.
For a prescribed macroscopic shear strain~$\overline{\gamma}$, the periodic fluctuating displacement field~$v(x)$ satisfies the equilibrium equation
\begin{equation} \label{eq:1d_bvp_coherent}
    \frac{\D}{\D x}\!\left(\mu(x)\frac{\D v}{\D x}\right) + \overline{\gamma}\,\frac{\D \mu(x)}{\D x} = 0 \, , \qquad x \in (0, \, L) \, ,
\end{equation}
cf.~\eqref{eq:homogBVP}.
The shear modulus is chosen as
\begin{equation} \label{eq:1d_mu_coherent}
    \mu(x) = \frac{\mu_0} {\alpha+(\alpha^{-1}-\alpha)\sin^2(\pi x/L)} \, ,
\end{equation}
where $\mu_0>0$ and $\alpha\in(0, \, 1)$ are material parameters.
Integrating~\eqref{eq:1d_bvp_coherent} once yields
\begin{equation}
    \frac{\D v}{\D x} = \frac{1}{\mu(x)} \left( C - \overline{\gamma}\mu(x) \right) \, ,
\end{equation}
where the integration constant~$C$ is determined from the periodicity condition on the fluctuating strain field, cf.~\eqref{eq:dispZeroAverage}, resulting in
\begin{equation}
    C = \frac{2\overline{\gamma}\mu_0\alpha}{1+\alpha^2} \, .
\end{equation}
The strain~$\gamma(x)$ is composed of the prescribed macroscopic strain~$\overline \gamma$ and the fluctuating strain, see~\eqref{eq:dispDecomposition}, and is given by 
\begin{equation} \label{eq:1d_gamma_exact}
    \gamma(x) = \frac{C}{\mu(x)} = \frac{2\overline{\gamma}\mu_0\alpha}{(1+\alpha^2)\mu(x)} \, .
\end{equation}
For the numerical experiments in the following two subsections, we choose the RVE parameters as
\begin{equation*}
    L=1, \, \mu_0 = 1,\text{ and } \, \alpha = \frac{3}{4} \, .
\end{equation*}
The respective strain field is given by
\begin{equation}
	\gamma(x) = \overline \gamma \left ( 1 - \frac{7}{25} \cos (2 \pi x) \right )  \, .
\end{equation}
%

%
%--------------------------------------------------------------------------------
\subsubsection{Single macroscopic FE strain \label{sec:example_1d_benchmark}}
\begin{figure}[bt]
        \subfloat[][Shear modulus $\mu(x)$ \label{fig:1d-figure_a}]{
        		\includegraphics[height=0.35\textwidth, keepaspectratio]{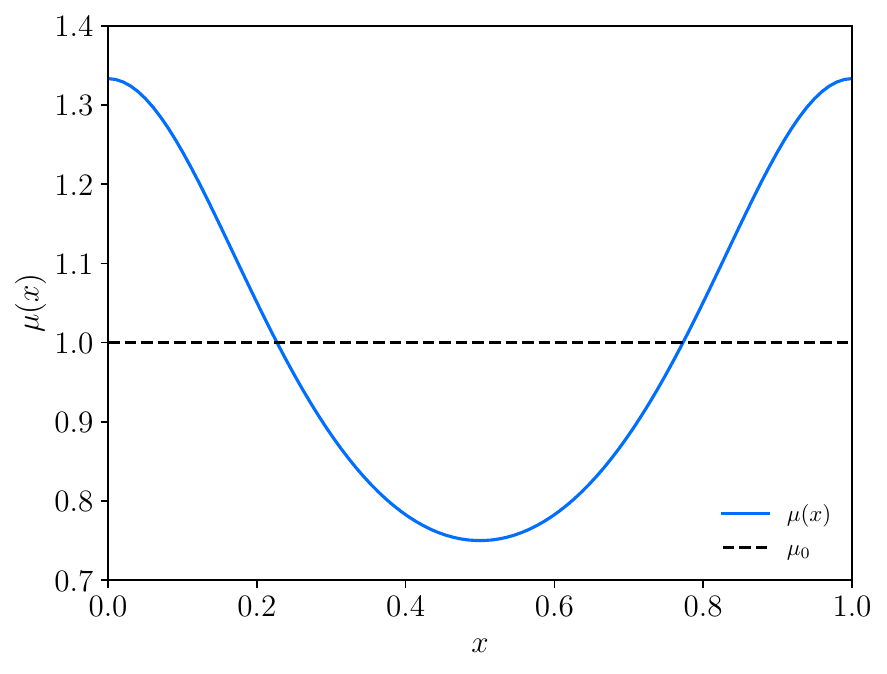}
        }
        \hfill
        \subfloat[][Strain $\gamma(x)$ \label{fig:1d-figure_b}]{
        		\includegraphics[height=0.35\textwidth, keepaspectratio]{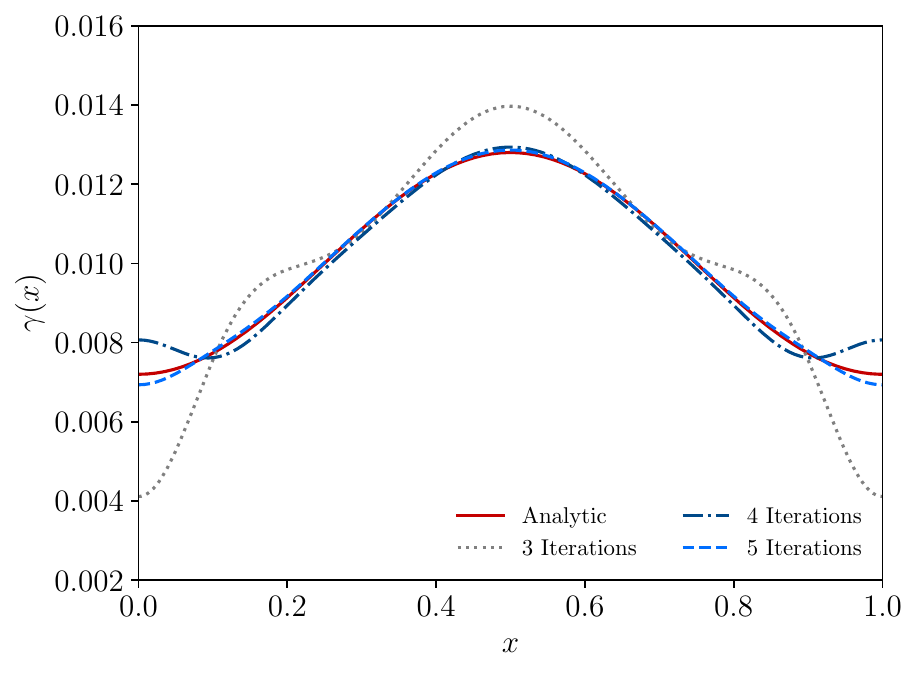}
        }
	\\ 
    	\vspace{0.4cm}
        \subfloat[][Relative $L_2$-norm error in strain~$\gamma(x)$ \label{fig:1d-figure_c}]{
        		\includegraphics[height=0.35\textwidth, keepaspectratio]{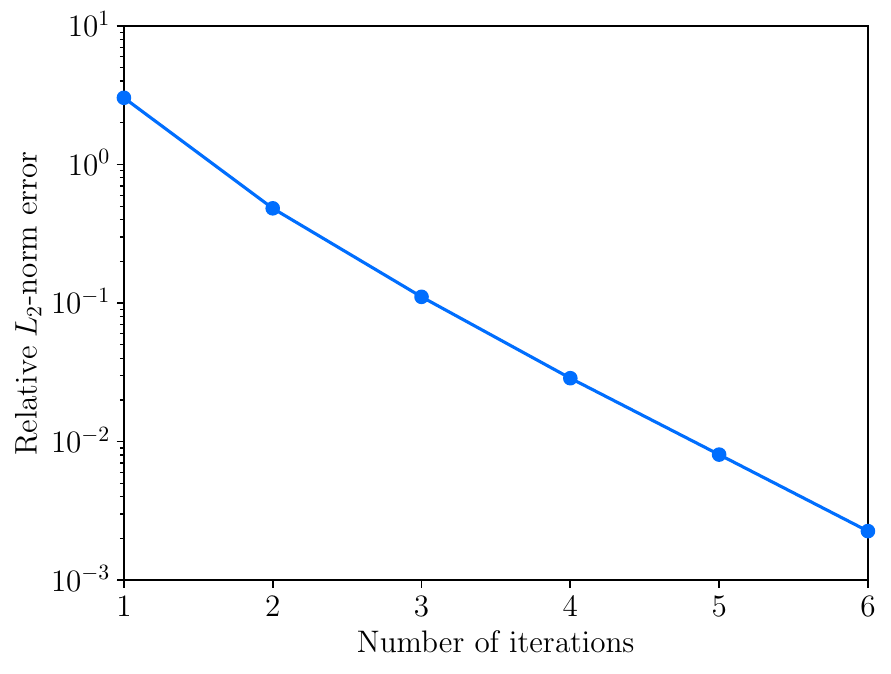}
        }
    	\hfill
        \subfloat[][Number of universal gates \label{fig:1d-figure_d}]{
        		\includegraphics[height=0.35\textwidth, keepaspectratio]{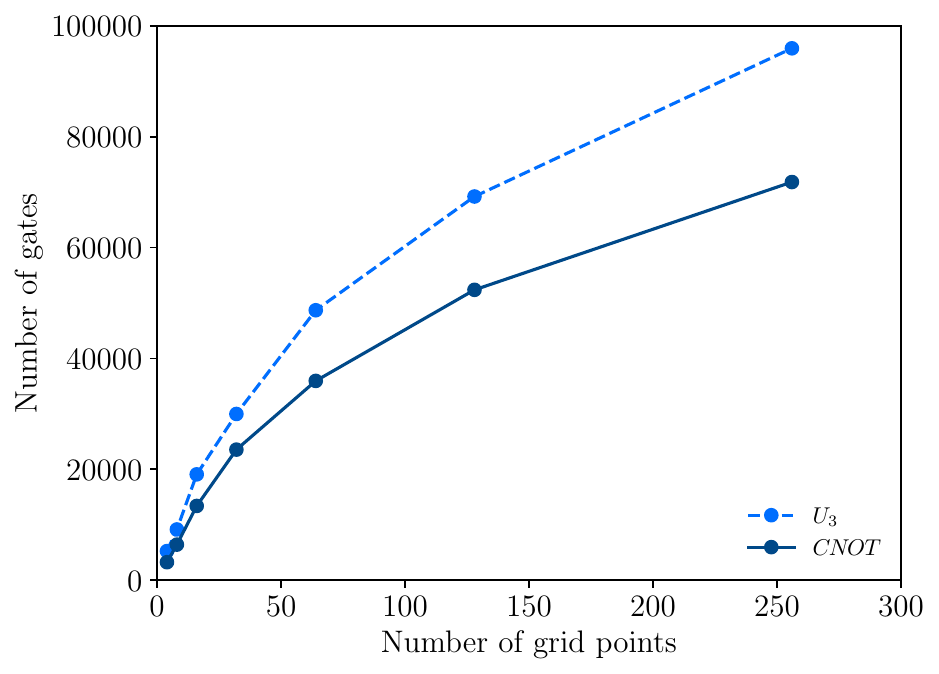}
        }
    \caption{One-dimensional RVE and a single prescribed macroscopic strain. (a) Spatial distribution of the shear modulus $\mu(x)$ and the reference shear modulus~$\mu_0$. (b) Exact and computed strains $\gamma(x)$ for iterations~$s \in \{3, \, 4, \, 5 \}$. (c) Convergence of the relative $L_2$-norm error in the computed strain.  (d) The total number of~$U_3$ and $\mathit{CNOT}$ gates.}
    \label{fig:1d-figure}
\end{figure}

%--------------------------------------------------------------------------------
%
We proceed to assess the accuracy and computational complexity of the proposed quantum approach for solving RVE problems in the case of a single prescribed macroscopic FE strain of $\overline {\gamma}=0.01$. In Figures~\ref{fig:1d-figure_a} and~\ref{fig:1d-figure_b}, the shear modulus and the exact strain field, respectively, are shown. For discretising the domain we consider uniform grids with~$N \in \{ 2^2, \, 2^3, \, 2^4 , \, 2^5, \, 2^6, \, 2^7, \, 2^8 \}$ grid points. The RVE problems are solved with the one-dimensional implementation of the quantum circuit introduced in Section~\ref{sec:incrementalUpdate}. The one-dimensional strain Green's function in the Fourier space is given by~$\hat \Gamma_k = -1 / \mu_0$. The modified shear modulus~$\mu (x) - \mu_0$ is first approximated on a classical computer by a polynomial of degree~$8$ via least-squares fitting and subsequently encoded into the quantum circuit following the approach of~\cite{liu2024towards}. The convergence of the computed strain field toward the exact strain field with increasing number of iterations for~$s \in \{3, \, 4, \, 5 \}$ is evident from Figure~\ref{fig:1d-figure_b}. The same can be deduced from the convergence of the relative $L_2$-norm error in the computed strain in Figure~\ref{fig:1d-figure_c}. The scaling of the total number of~$U_3$ and $\mathit{CNOT}$ gates with respect to the number of grid points~$N$ is depicted in Figure~\ref{fig:1d-figure_d}. The number of gates depends polylogarithmically on the number of grid points, consistent with the complexity $O(\log^c N)$ of the proposed quantum approach.

%
%--------------------------------------------------------------------------------
\subsubsection{Many macroscopic FE strains \label{sec:example_1d_multi}}
%--------------------------------------------------------------------------------
%
We now consider the solution of the introduced one-dimensional RVE problem for many, i.e.~$M \gg 1$,  prescribed macroscopic FE strains. The discretisation, the shear modulus and the number of fixed-point iteration steps are the same for all the RVE problems. According to the results of the previous section, the sequential solution of all the~$M$ RVE problems on a quantum computer will have a computational complexity~$O(M \log^c N)$. In the \emph{QAFE$^2$} framework, all RVE problems are solved simultaneously, resulting in substantially lower computational complexity. As introduced in Section~\ref{sec:quantum_parallel_rve_update},  in \emph{QAFE$^2$}, after encoding the prescribed macroscopic strains serially in the initialisation step, all the RVE problems are solved simultaneously by employing the standard quantum RVE solution algorithm. It bears emphasis that, for a reasonably fine discretisation, encoding the macroscopic strain is significantly less expensive than solving the RVE problem using fixed-point iteration.

\begin{figure}
    \centering
        \subfloat[][Number of grid points $N=2^4$ \label{fig:1d-ensemble-gatecount-n4}]{
        \includegraphics[height=0.35\textwidth, keepaspectratio]{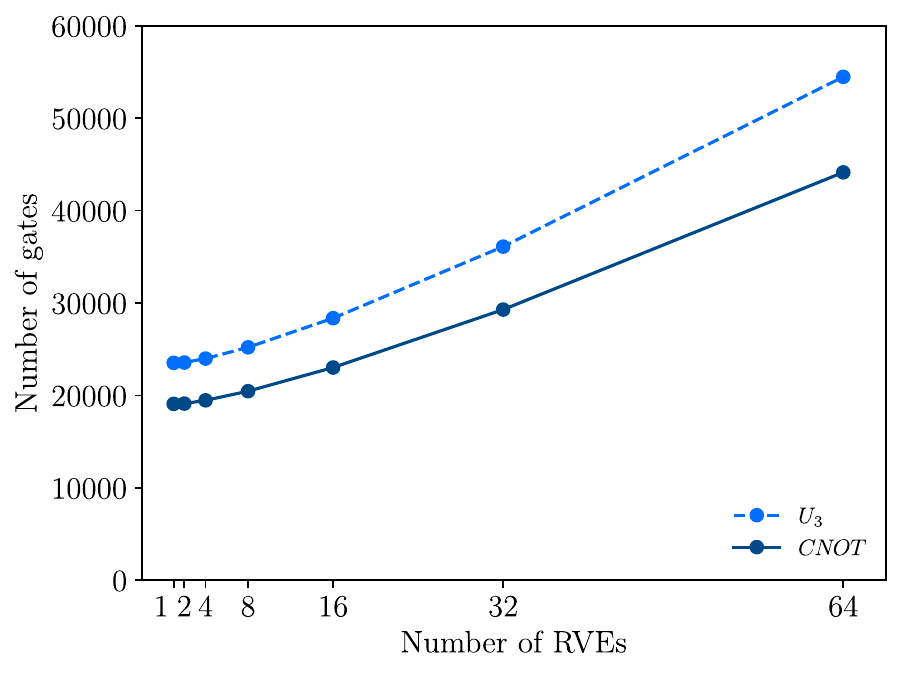}
        }
        \hfill
        \subfloat[][Number of grid points $N=2^{10}$  \label{fig:1d-ensemble-gatecount-n10}]{
        \includegraphics[height=0.35\textwidth, keepaspectratio]{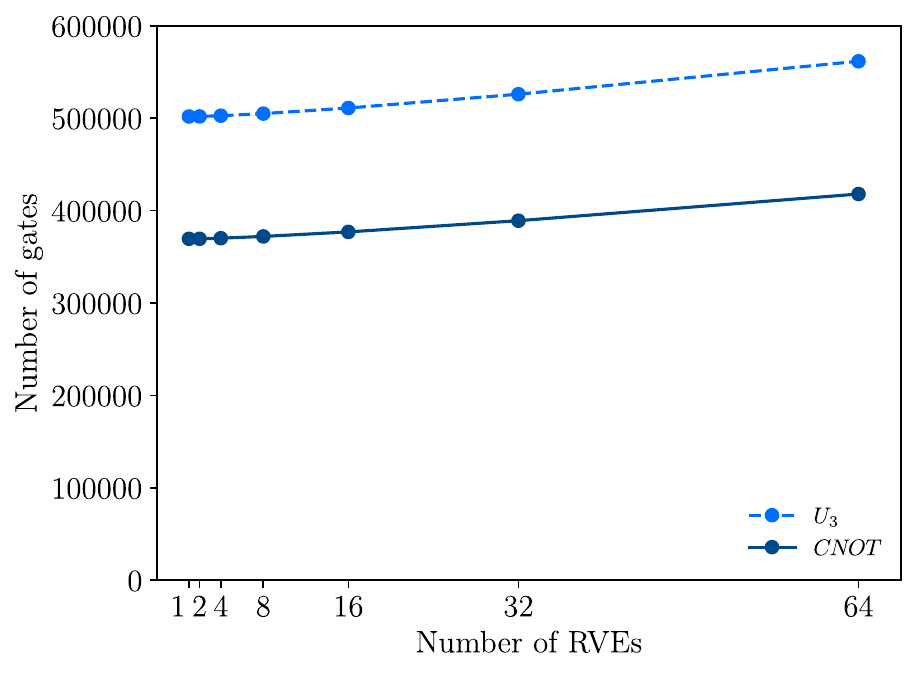}
        }
    \caption{One-dimensional RVE with~$M$ prescribed macroscopic strains. Total number of $U_3$ and $\mathit{CNOT}$ gates for two different discretisations with $N=2^4$ and $N=2^{10}$. \label{fig:1d-ensemble-gatecount}}
\end{figure}

To assess the computational complexity of solving the RVE problems in the \emph{QAFE$^2$} framework, we consider discretisations with a resolution of~$ N=2^4$ and~$N=2^{10}$ grid points.  In Figure~\ref{fig:1d-ensemble-gatecount}, the scaling of the number of~$U_3$ and $\mathit{CNOT}$ gates with respect to the number of RVE problems is plotted. Clearly, for both discretisations, the number of required gates is significantly smaller than solving $M$ different RVE problems sequentially. The two plots indicate a complexity of~$O (M\log^{c} M + \log^{c} N)$, which can also be straightforwardly theoretically shown. The multi-controlled gates for encoding the~$M$  macroscopic strains have a complexity~$O (M\log^{c} M)$, and the solution of the RVE problems has a complexity~$O (\log^{c} N)$. 
%
%--------------------------------------------------------------------------------
\subsection{Two-dimensional RVE \label{sec:example_2d_separable}}
\begin{figure}[]
        \subfloat[][Shear modulus field $\mu(\vec x)$ \label{fig:2d-mu}]{
        \includegraphics[width=0.32\textwidth]{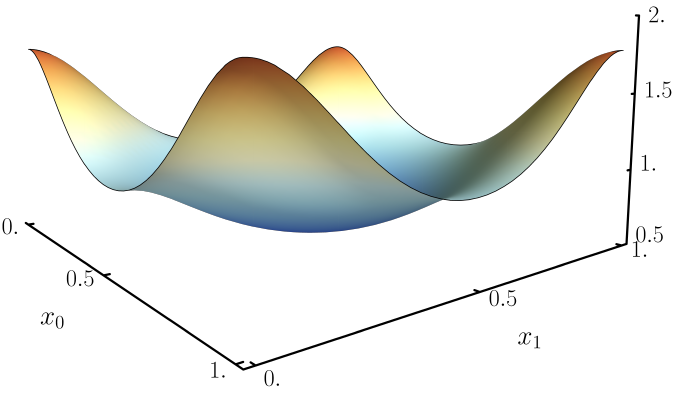}
        }
    \hfill
        \subfloat[][Exact strain $\gamma_0(\vec x)$ \label{fig:2d-gamma0-exact}]{
        \includegraphics[width=0.32\textwidth]{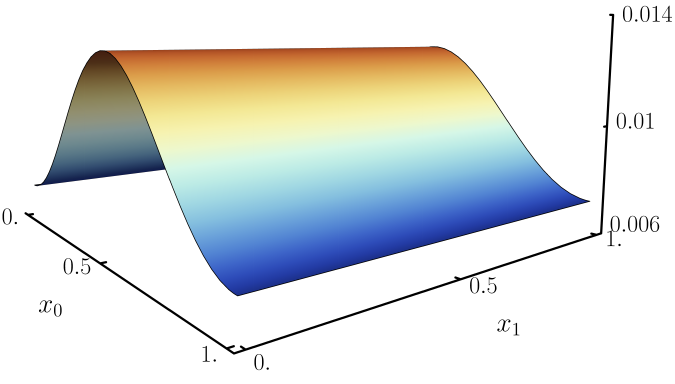}
        }
    \hfill
        \subfloat[][Exact strain $\gamma_1(\vec x)$  \label{fig:2d-gamma1-exact}]{
        \includegraphics[width=0.32\textwidth]{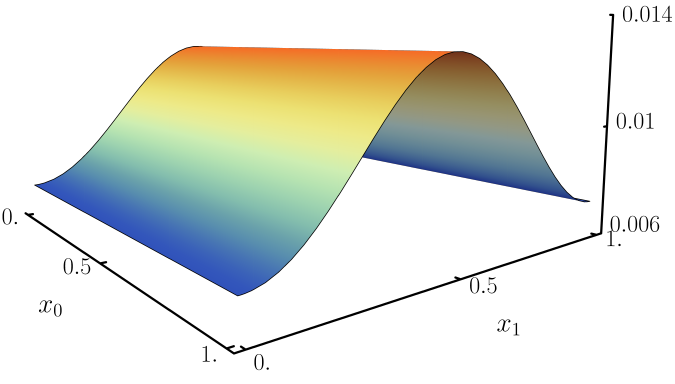}
        }
    \caption{Two-dimensional RVE. (a) Chosen shear modulus~$\mu(\vec x)$. (b, c) Strain components~$\gamma_0(\vec x)$ and $\gamma_1(\vec x)$ of the analytically obtained exact strain field.}
    \label{fig:2d-analytical}
\end{figure}

%--------------------------------------------------------------------------------
%
Our next example concerns a two-dimensional RVE defined on $ \Omega=(0, \, L) \times (0, \, L) $ with periodic boundary conditions. The spatial coordinates are denoted by \mbox{$\vec x = ( x_0 \quad x_1 )^\trans$}, and the periodic fluctuating displacement field \mbox{$v(\vec x)= (  v_0(\vec x) \quad v_1(\vec x) )^\trans$} satisfies the equilibrium equation~\eqref{eq:homogBVP}. We choose a periodic shear modulus~$\mu(\vec x)$, which is separable and is given by
\begin{subequations}
\begin{align} \label{eq:2d_mu_separable}
	\mu(\vec x) & =\kappa_0(x_0) \kappa_1 (x_1) \, ,
\intertext{where}  
	\kappa_0(x_0) &=
	\frac{\mu_0} {\alpha+(\alpha^{-1}-\alpha)\sin^2(\pi x_0/L)} \, ,  \\
	\kappa_1(x_1) &=
	\frac{\mu_0} {\alpha+(\alpha^{-1}-\alpha)\sin^2(\pi x_1/L)} \, .
\end{align}
\end{subequations}
The resulting boundary value problem~\eqref{eq:homogBVP} is fully heterogeneous while remaining analytically tractable. For a separable shear modulus, the fluctuating displacement field admits the additive representation
\begin{equation} \label{eq:fluct_disp_additive}
	v(\vec x)= w_0(x_0) + w_1(x_1) +C \, .
\end{equation}
with a constant~$C$ fixed by the zero-mean condition.
For the computations, we choose the RVE parameters and the prescribed macroscopic strain as 
\begin{equation*}
L=1, \, \mu_0 = 1, \, \alpha = \frac{3}{4}, \text{ and } \,  \overline {\vec \gamma} = \begin{pmatrix}0.01  & 0.01 \end{pmatrix}^\trans \, .
\end{equation*}
Using the ansatz~\eqref{eq:fluct_disp_additive}, the corresponding strain field can be derived in closed form as
\begin{equation}
	\vec \gamma (\vec x) =
	\begin{pmatrix}
    \dfrac{25 - 7\cos(2\pi x_0)}{2500}\\[6pt]
    \dfrac{25 - 7\cos(2\pi x_1)}{2500}
    \end{pmatrix} \, .
\end{equation}
The chosen shear modulus field $\mu(\vec x)$ and the corresponding analytical strain components $\gamma_0(\vec x)$ and $\gamma_1(\vec x)$ are visualised in Figure~\ref{fig:2d-analytical}. Although the analytical strain field admits a separable representation, the quantum solver operates on the full two-dimensional problem, and all intermediate strain fields generated during the fixed-point iteration remain genuinely two-dimensional. This separability provides a convenient internal consistency check, enabling the convergence of each strain component to be assessed independently along one-dimensional slices.

\begin{figure}[tb]
        \subfloat[][Strain $\gamma_0(x_0, x_1=0.5)$ \label{fig:2d-slices-0}]{
        \includegraphics[height=0.35\textwidth, keepaspectratio]{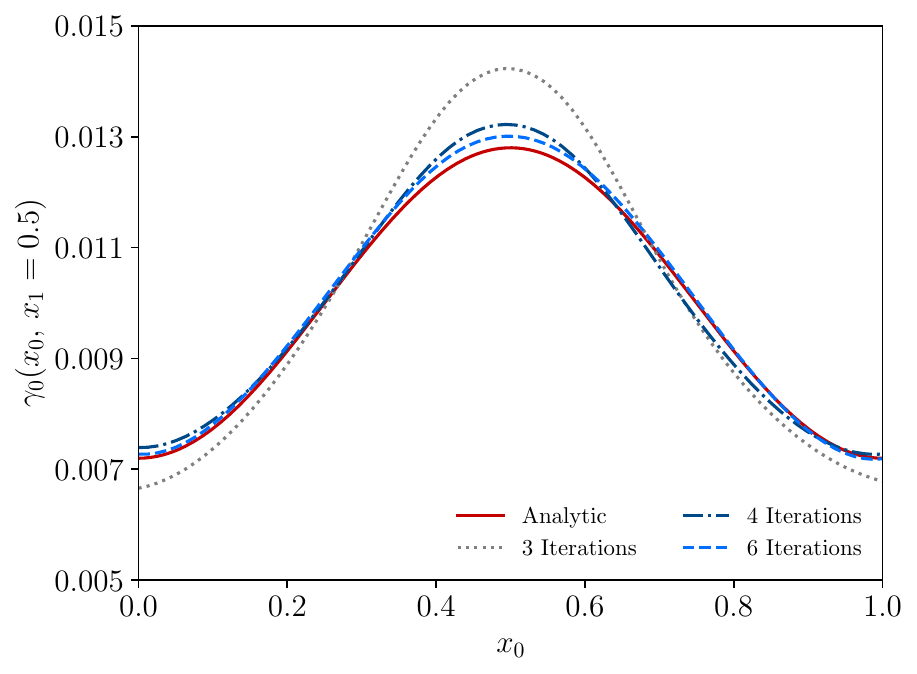}
        }    
        \hfill
        \subfloat[][Strain $\gamma_1(x_0, x_1=0.5)$  \label{fig:2d-slices-1}]{
        \includegraphics[height=0.35\textwidth, keepaspectratio]{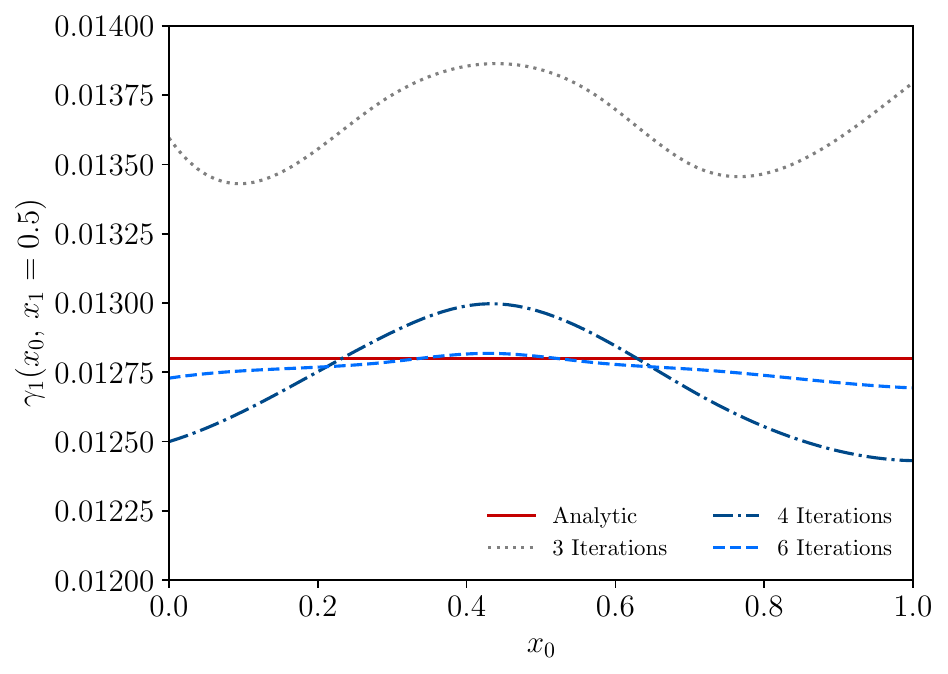}
        }
        \\ 
        \subfloat[][Relative $L_2$-norm error in strain~$\vec \gamma (\vec x)$    \label{fig:2d-residual}]{
        \includegraphics[height=0.35\textwidth, keepaspectratio]{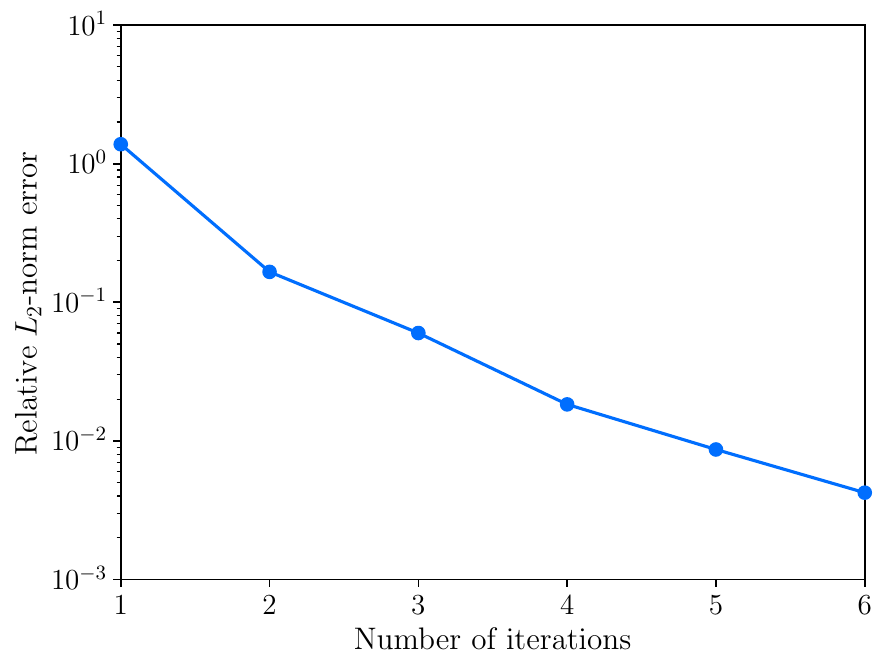}
        }
    \hfill
        \subfloat[][Number of universal gates \label{fig:2d-gatecount}]{
        \includegraphics[height=0.35\textwidth, keepaspectratio]{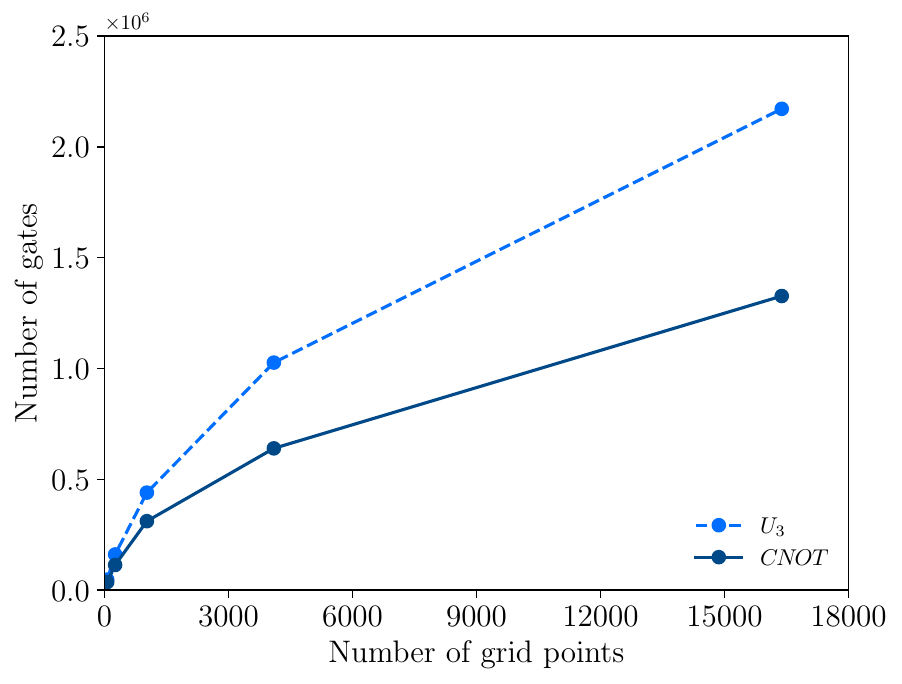}
        }
    \caption{Two-dimensional RVE and a single prescribed macroscopic strain. (a, b) Exact and computed strain components $\gamma_0 (x_0, x_1=0.5)$ and~$\gamma_1(x_0, x_1=0.5)$ along the horizontal at~$x_1 =0.5$ for iterations~$s \in \{3, \, 4, \, 6 \}$. (c) Convergence of the relative $L_2$-norm error in the computed strain. (d) Total number of~$U_3$ and $\mathit{CNOT}$ gates. \label{fig:2d-figure}}
\end{figure}

For discretising the domain we use grids with $N\times N$ grid points, where ~$N \in \{ 2^2, \, 2^3, \, 2^4 , \, 2^5, \, 2^6, \, 2^7\}$. The RVE problems are solved using the quantum circuit introduced in Section~\ref{sec:incrementalUpdate}. As discussed in Section~\ref{sec:incrementalFourierSpace}, the two-dimensional strain Green's operator in Fourier space $\hat{\Gamma}_{k^0 k^1}$ is block-encoded via the LCU approach. Bivariate polynomials of degrees at most~$7$ and~$6$ are used to approximate the scalar coefficient functions~$ \alpha_1 (k^0, k^1)$ and~$ \alpha_2 (k^0, k^1)$ associated with the Pauli~$X$ and~$Z$ gates, cf.~\eqref{eq:gammaLCU_0}. These polynomials are then quantum encoded following the approach of~\cite{liu2024towards}. When the extended-domain method described in~\ref{app:extended-domain} is employed, these polynomial degrees are reduced to~$3$ and~$4$. The spatial fluctuation of the shear modulus relative to its reference values is encoded into the quantum circuit using a bivariate polynomial approximation of degree at most~$4$ in each spatial coordinate. All polynomial coefficients are obtained via least-squares fitting on a classical computer. 
The convergence of the computed strain towards the exact strain along the horizontal midline for iteration steps~$s \in \{ 3, \, 4, \, 6\}$ is shown in Figures~\ref{fig:2d-slices-0} and~\ref{fig:2d-slices-1}.  The convergence of the relative $L_2$-norm error in the computed strain with increasing number of iterations is shown in Figure~\ref{fig:2d-residual}. The scaling of the total number of~$U_3$ and $\mathit{CNOT}$ gates with respect to the number of grid points~$N^2$ is depicted in Figure~\ref{fig:2d-gatecount}. The number of gates depends polylogarithmically on the number of grid points, consistent with the polylogarithmic complexity~$O(\log^c N)$ of the proposed quantum approach. 

%
%--------------------------------------------------------------------------------
\section{Conclusions \label{sec:conclusions}}
%--------------------------------------------------------------------------------
%
We have introduced a quantum--classical framework, termed \emph{QAFE$^2$}, for concurrent multiscale finite element analysis that fundamentally alters the computational scaling of RVE-based homogenisation. At the single-RVE level, the approach builds on the classical FFT-based fixed-point scheme of Moulinec and Suquet~\cite{moulinec1998} and provides a complete quantum reformulation in which the core operations---Fourier transformation, constitutive updates, and fixed-point iteration---are implemented by quantum circuits. By exploiting the polylogarithmic complexity of the QFT and block-encoded representations of the strain Green's operator, the cost of solving an individual RVE discretised on an $N\times N$ grid is reduced from $\mathcal{O}(N^2\log N)$ to $\mathcal{O}(\log^c N)$ for a modest constant $c$, representing an \emph{exponential speedup} relative to the best available classical algorithms.

Beyond the acceleration of single-RVE problems, the defining contribution of \emph{QAFE$^2$} lies in its treatment of \emph{concurrency}. By encoding the macroscopic strains associated with all quadrature points into orthogonal subspaces of a single quantum state, the framework enables the simultaneous solution of all $M$ RVEs in a finite element problem using a single instance of the fixed-point iteration circuit. The resulting overall complexity scales essentially linearly with $M$ only through the state-initialisation overhead, while the microscopic solve itself is performed once, independently of $M$. This form of quantum parallelisation has no classical analogue and directly targets the dominant bottleneck that has historically limited the practical applicability of FE$^2$ methods in large-scale simulations.

The presented numerical experiments for one- and two-dimensional RVEs for a linear model problem in antiplane shear verify both the accuracy of the proposed quantum fixed-point iteration and the predicted scaling behaviour. For problems with known analytical solutions, the quantum formulation reproduces the exact homogenised response while exhibiting gate-count growth consistent with theoretical complexity estimates. Although the simulations are carried out on noiseless quantum emulators, they demonstrate that the algorithmic building blocks of \emph{QAFE$^2$} can be assembled into concrete, end-to-end circuits whose depth and ancilla requirements remain independent of the microscopic resolution.

Notwithstanding these advances, several notable limitations and open challenges remain. The measurement of the computed average microscopic stress by repeated circuit evaluations requires post-selection of the physically relevant states. The number of discarded states increases with the number of iterations and RVEs, leading to reduced success probabilities. This may require additional techniques, such as amplitude amplification~\cite{brassard2002quantum}, to improve the efficiency of the measurement procedure. In addition, our analysis assumes ideal, fault-tolerant quantum hardware and periodic RVEs with relatively simple constitutive behaviour. Extending the framework to more general boundary conditions, nonlinear material models, and inelastic or history-dependent responses will require further developments in quantum encoding and circuit design. Likewise, the impact of noise, finite coherence times, and error correction on the proposed algorithms remains to be quantified.

Notwithstanding these challenges, the results to date suggest that quantum computing offers more than incremental acceleration for computational mechanics on future fault-tolerant quantum computers. By reshaping the algorithmic structure of multiscale analysis and enabling genuinely concurrent microscopic solves, \emph{QAFE$^2$} points towards a new paradigm in which quantum hardware acts as a constitutive engine embedded within otherwise classical finite element workflows. Recent progress in quantum error correction and logical qubits provides a clear motivation for developing such algorithms already now~\cite{bluvstein2024logical,google2025quantum}. As quantum technologies mature, hybrid approaches may ultimately enable fully concurrent multiscale simulations at scales far beyond the reach of classical computing.

\section*{Acknowledgments}

M.O. gratefully acknowledges the financial support of the Centre Internacional de M\`etodes Num\`erics a l’Enginyeria (CIMNE) of the Universitat Polit\`ecnica de Catalunya (UPC), Spain, through the UNESCO Chair in Numerical Methods in Engineering. 

\appendix
%
%--------------------------------------------------------------------------------
\section{Encoding of piecewise discontinuous functions \label{app:extended-domain}}
\begin{figure}
    \centering
    \includegraphics[scale=0.9]{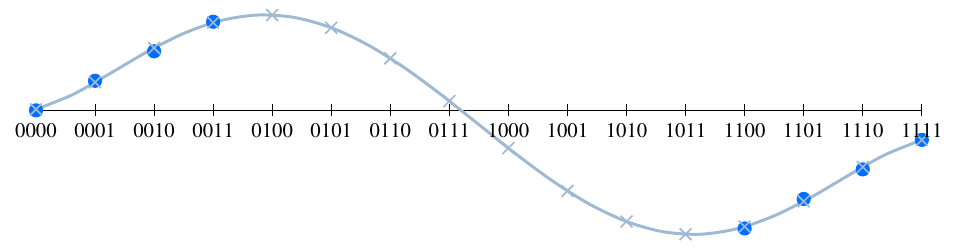} 
     \caption{Approximation of the function~$\tilde r(k)$ defined in the intervals~$[0, \, 4)$ and~$[12, \, 15]$ using a quintic polynomial. The dots denote the grid-point values of~$\tilde{r}(k)$, and the crosses the grid-point values of the quintic polynomial (solid line), obtained by least-squares fitting. The grid points are indexed in binary.}
\label{fig:discontinuous-extended-before}
\end{figure}

%--------------------------------------------------------------------------------
%
Piecewise discontinuous functions, like the bivariate LCU coefficient~$\alpha_1 (k^0, k^1)$ in~\eqref{eq:gammaLCU_0},  can be efficiently quantum-encoded using a single polynomial by embedding the problem in an extended domain. We illustrate the construction using the univariate relabelling function~$r(k)$ introduced in~\eqref{eq:componentsToFrequencies}  on a grid with~$N=8$ grid points. Consider the extended domain with~$N=16$ grid points and the auxiliary function
\begin{equation} \label{eq:componentsToFrequenciesTilde}
	\tilde{r}(k) = \begin{cases}
		k & 0 \le k < 4\\
		k - 16 & 12 \le k \le 15 \, .
	\end{cases}
\end{equation}
The values of $\tilde{r}(k) $ between~$ 4 \le k  < 12$ are unspecified.  In Figure~\ref{fig:discontinuous-extended-before}, the values of~$\tilde r(k)$ at grid points, and its least-squares approximation using a polynomial of degree five are shown. The polynomial approximant is continuous and infinitely smooth. We seek to determine from the polynomial approximant an approximant for~$r(k)$ defined in the interval~$[0, \, 7] $, which has a discontinuity at~$k=4$. To this end, it is expedient to write the grid indices in binary as in Figure~\ref{fig:discontinuous-extended-before}. By inspection, flipping the value of the leftmost bit when the second bit from the left is in state~$\ket 1$, yields an approximant for~$r(k)$ in the interval~$[0, \, 8)$ with a discontinuity at~$k=4$.

\begin{figure}[b]
\centering
\begin{minipage}[b]{0.375\textwidth}
\centering
\subfloat[][Encoding of~$r(k)$   \label{fig:extended-domain-circ}]{
	\includegraphics[scale=0.9]{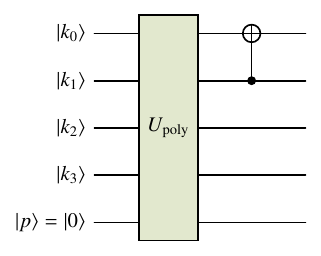} 
}
\end{minipage}
\hspace{0.05\textwidth}
\begin{minipage}[b]{0.475\textwidth}
\centering
\subfloat[][$U_{\text{exch}}$ mapping~$\alpha \ket{000} + \beta \ket{111}$ to~$\beta \ket{000} + \alpha \ket{111}$ \label{fig:gray-code-exch}]{
	\includegraphics[scale=0.9]{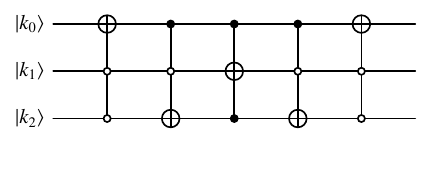} 
}
 \end{minipage}
\caption{Auxiliary quantum circuits.}
\end{figure}

A quantum circuit implementation of the sketched construction is shown in Figure~\ref{fig:extended-domain-circ}. The unitary $U_{\text{poly}}$ encodes the polynomial of degree five that has been classically determined by least-squares fitting to~$\tilde r(k)$.   The~$\mathit{CNOT}$ gate flips the state of the qubit~$\ket {k_0}$ when the qubit~$\ket{k_1}$ is in state~$\ket {1}$. The values of~$r(k)$  correspond to the states with~$\ket{0} \ket{k_1 k_2 k_3 } \ket{1}$. For possible implementations of~$U_{\text{poly}}$,  see~\cite{liu2024towards,febrianto2025quantum,rosenkranz2025quantum}. 

%--------------------------------------------------------------------------------
\section{Component exchange unitary  \label{app:graycode_exch}}
%--------------------------------------------------------------------------------
%
The unitary~$U_{\text{exch}}$ introduced in~\eqref{eq:exch_bar_gamma} implements the permutation
\begin{equation}
U_{\text{exch}}  \colon \alpha \ket{k_a} + \beta \ket{k_b} \mapsto  \beta \ket{k_a} + \alpha \ket{k_b} \, ,
\end{equation}
where~$k_a$ and~$k_b$ are two components of the computational basis and~$\alpha$ and~$\beta$ their coefficients. All other components of the state vector remain unaffected. Hence, the unitary~$U_{\text{exch}} $ acts as an~$X$ gate in the plane spanned by the basis vectors~$\ket{k_a}$ and~$\ket{k_b}$. For arbitrary~$\ket{k_a}$ and~$\ket {k_b}$, this operation can be implemented using multi-controlled $X$ gates based on a Gray code connecting~$k_a$ to~$k_b$~\cite{ikeAndMike}. A Gray code is a sequence of binary numbers in which successive members differ exactly in a single bit. For instance, for~$k_a = 000$  and~$k_b=111$ a possible Gray code sequence is 
\begin{equation}
	000, 100, 101, 111 \, .
\end{equation}
The unitary~$U_{\text{exch}}$ can be implemented with five multi-controlled $X$ gates. Three of the gates implement the sequence~$000, 100, 101, 111$ and the other two the reverse sequence~$101, 100, 000$. As shown in Figure~\ref{fig:gray-code-exch}, at each step an~$X$  gate is applied to the differing qubit conditioned on all other qubits being in the same state. See~\cite{ikeAndMike} for further details.

\bibliographystyle{elsarticle-num-names}
\bibliography{quantumAccelerated}

\end{document}